\documentclass{amsart}
\usepackage{amssymb}
\usepackage{amscd}
\DeclareFontFamily{OMS}{cmsy}{%
\fontdimen16\font=3pt
\fontdimen17\font=3pt}
\makeatletter
\renewcommand{\subsection}{\@startsection{subsection}{2}{\z@}%
{\baselineskip}{0.5\baselineskip}{\bfseries}}
\def\l@section{\def\@tocpagenum##1{\hss{\bfseries ##1}}%
\@tocline{1}{8pt}{0pc}{}{\bfseries}}
\def\l@subsection{\def\@tocpagenum##1{}
\@tocline{2}{2pt}{2pc}{2pc}{}}
\makeatother
\def\dj{d\kern-.30em\raise1.25ex\vbox{\hrule width .3em height .03em}}
\def\Dj{D\rlap{\kern-.70em\raise0.75ex
\vbox{\hrule width .3em height .03em}}}
\def\dM{\mbox{$\smash{\vphantom{d}^{\raise0.4ex\hbox{$\scriptscriptstyle M$}}%
\!\!\!d}$}}
\def\bla#1{$(${\it #1\/{}}$)$}
\def\U{\mathrm{U}} 
\def\O{\mathrm{O}} 
\def\SU{\mathrm{SU}}
\def\SO{\mathrm{SO}}
\def\restr{\restriction}

\def\cal{\mathcal}
\def\der{\mathfrak{der}}
\def\Bbb{\mathbb}
\def\e{\epsilon}
\def\k{\kappa}
\def\rig{\wp}
\def\tr{\mathrm{tr}}
\def\ad{\mathrm{ad}}
\def\id{\mathrm{id}}

\def\zh{\mathfrak{zh}}

\def\hor{\mathfrak{hor}}
\def\V{\Bbb{V}}
\def\v{\varkappa}
\def\Fh{F^\wedge}
\def\inv{{i\!\hspace{0.8pt}n\!\hspace{0.6pt}v}}
\def\im{\mathrm{im}}
\def\cstV{\widehat{\cal{V}}}
\def\cstB{\widehat{\cal{B}}}
\def\vint{\int_{\!\uparrow}}
\def\modM{\Theta} 
\def\modP{\widehat{\Theta}}
\def\modH{\Lambda}
\def\intP{\int_{\!P}}
\def\hstr{{\star}} 
\def\hstrP{\hstr_P}
\def\hstrM{\hstr_M}
\def\ll{{\lozenge}}

\def\lin{\mathrm{lin}}
\def\D{\mathbb{D}}

\def\LapP{\Delta_P}
\def\LapM{\Delta_M}
\def\LapS{\Delta_{\cal{S}}}
\def\ClP{\mathfrak{cl}[P]}
\def\ClM{\mathfrak{cl}[M]}
\def\Fcl{F_{\mathrm{cliff}}}
\def\Fspi{F_{\cal{S}}}
\def\Salg{\mathbb{S}}
\def\vspi{\v_\Salg}
\def\vcl{\v_*}
\def\SpiM{\cal{S}_M}
\def\ClV{\mathrm{cl}[\V,g,\sigma,\mc]}
\def\map#1{#1_\star}
\def\amaph#1{*_{#1}}
\def\Mor{\mathrm{Mor}}
\def\RepG{R(G)}
\def\hbim#1{\cal{F}_{#1}}
\def\bim#1{\cal{E}_{#1}}
\def\Sum{{\displaystyle\sum}}
\def\mc{\Sigma}
\def\mcVt{\nu_\mc}
\def\Vmc{\V_{\mc}}
\def\horPmc{\hor_{P,\mc}}
\def\Wmc{\Omega_{M,\mc}}
\def\vmc{\v_{\mc}}
\renewcommand{\thepage}{\ifnum\value{page}=1 \else\arabic{page}\fi}
\newtheorem{thm}{Theorem}[section]
\newtheorem{pro}[thm]{Proposition}
\newtheorem{lem}[thm]{Lemma}
\theoremstyle{definition}
\newtheorem{defn}{Definition}
\newenvironment{pf}{\proof[\proofname]}{\endproof}
\numberwithin{equation}{section}
\thanks{Supported by the United States-Mexico Foundation for Science. Revision 7, April 2000}
\begin{document}

\title[quantum principal bundles]{General Spinor Structures\\On Quantum Spaces}
\author{micho {\Dj}UR{\Dj}EVICH}
\address{Instituto de Matematicas,
UNAM, Area de la Investigacion Cientifica,
Circuito Exterior, Ciudad Universitaria, M\'exico DF, cp 04510,
MEXICO} 
\email{ micho@matem.unam.mx\newline
\indent\indent http://www.matem.unam.mx/{\~\/}micho}
\begin{abstract}
A general theory of quantum spinor structures on quantum spaces is presented, 
within the conceptual framework of the formalism of quantum principal bundles. 
Quantum analogs of all basic objects of the classical theory are constructed
and analyzed. This includes Laplace and Dirac operators, quantum 
versions of Clifford and spinor bundles, a Hodge $*$-operator, appropriate 
integration operators, and mutual relations of these objects.
We also present a self-contained formalism of braided Clifford algebras. 
Quantum phenomena appearing in the theory are discussed, including a very 
interesting example of the Dirac operator associated to a quantum Hopf fibration. 
\end{abstract}
\maketitle
\enlargethispage{2\baselineskip}
\tableofcontents
\filbreak
\section{Introduction}

Classical theory of spinors is a cornerstone of important bridges that
interconnect worlds of geometry, algebra and physics. 

Just as an illustration---let us recall that using spinor structures, 
it is possible to address various fundamental topics of Riemannian 
geometry in an elegant and effective way. The Dirac operator is intrinsically 
related with a general index theory of elliptic operators. Spin bundles 
provide powerful tools for the study of gauge field theories with fermions.  

Classical geometry is a very special case of quantum 
($\Leftrightarrow$ non-commutative) geometry, which
incorporates many ideas of quantum physics into the world of geometry. 
And if we try to understand the structure of physical space-time at 
very small scales, it becomes clear that classical concepts of the 
space-time continuum loose their validity. 

It is natural to look for a quantum version of the spinor theory, hoping that 
it would be at least as interesting as its classical counterpart. 

The aim of this study is to present a general theory of spinor structures
over quantum spaces, in the spirit of non-commutative
differential geometry \cite{c}. The main conceptual framework for 
our considerations is the theory of quantum principal bundles \cite{d1,d2},
where quantum groups play the role of structure groups and 
general quantum spaces play the role of base manifolds. 

The formalism presented here could be used as a possible tool for 
developing a theory of fermions over a quantum space-time, which would
be appropriate at the level of ultra-small distances characterized by 
the Planck length. Our formulation fulfills various conditions proposed 
in a general axiomatic framework \cite{c}. However, the formulation of \cite{c} 
does not require quantum principal bundles as the basic underlying structure. 
Further, some key conditions of \cite{c} are broken  in our formalism. 
This includes the spectral asymptotics of the quantum Dirac operator, 
which in our case could be very different from the classical behavior. 
In particular, our constructions are not compatible with the Dixmier trace.  

As far as pure geometry is concerned, our constructions provide a 
coherent framework for the formulation of quantum elliptic operators 
\cite{d-elliptic} and the study of the corresponding index theorems. 

The results of this paper include as a special case the formalism of 
spin structures studied in \cite{d-clspi}, where quantum spin bundles with 
classical structure groups were considered (and it was assumed that
the differential calculus over the structure group is classical, too). 

The paper is organized as follows. 

\smallskip
In the next section we are going to introduce basic structural elements of
quantum Riemannian geometry. We shall use the general theory of frame structures
on quantum principal bundles \cite{d-frm2}, in order to develop the idea of a quantum 
space equipped with a metric. In particular, we shall explain how to construct a
graded *-algebra $\hor_P$ representing horizontal forms,  starting from `abstract 
coordinate 1-forms' and a quantum principal bundle $P$. In accordance with 
\cite{d-frm2}, the space of abstract coordinate 1-forms $\V$ will be defined
as the left-invariant part of a bicovariant bimodule $\Psi$ over the structure 
quantum group $G$. The group $G$ acts on $\V$ by `orthogonal transformations'. 
The space $\V$ carries a very interesting geometrical structure, 
and in particular \cite{w3} there exists a canonical braid operator 
$\tau\colon\V\otimes\V\rightarrow \V\otimes\V$ playing the role of the transposition 
map. Furthermore, we shall introduce abstract Levi-Civita connections. In a certain sense, 
these objects contain the whole geometrical information about quantum frame structures. 

This is the most subtle part of the formalism, as it requires
to introduce carefully a number of non-trivial conditions on the base space $M$, the structure quantum group $G$ and the bundle $P$. 

After presenting the main algebraic setup, and reviewing basic ideas of \cite{d-frm2}, 
we shall consider certain special conditions which will further justify
our geometrical interpretation of $M$---as a quantum space equipped with a metric. 

This includes
analytic conditions---the existence of the appropriate C*-algebraic completions of both 
the base space and the bundle *-algebras, as well as the existence of a suitable 
`homogeneous' measure on the base space $M$. Combining this measure with the integration
along the fibers of $P$ we shall construct a natural measure on the bundle. 

One of the main purely algebraic extra 
conditions will be the existence of a `volume element' in the algebra of coordinate 
horizontal forms. With the help of the volume element and the measure on $P$, 
it will be possible to construct the integration map 
$\intP\colon\hor_P\rightarrow\Bbb{C}$. 

A quantum version of the Euclidean structure on $\V$ will be represented by a
metric form $g\colon\V\otimes\V\rightarrow \mc$, where $\mc$ is a *-algebra
of {\it abstract metric tensor coefficients}. In general $\mc\neq\Bbb{C}$, and there
exist deep reasons why it is necessary to assume that
components of the metric generate a non-commutative *-algebra. 

After discussing the integration maps, we shall present the construction of the quantum Hodge $*$-operator, 
with the help of which it will be possible to express elegantly a canonical scalar product
of horizontal forms, along the lines of classical geometry. 
Next, we shall introduce the coderivative operator, the Laplace operator and analyze their  
properties. Explicit coordinate formulas for the Laplacian and the adjoint derivative will be 
derived, too. 

In Section~3 we shall introduce the concept of a quantum spinor structure, in a complete
analogy with classical geometry. These structures will be defined as certain `covering bundles' 
of the `true orthonormal frame bundles'. Their structure group
will be a kind of a quantum spin group. We shall proceed by introducing quantum versions 
of the Clifford bundle algebra and the associated spinor bundle.

Section~4 is devoted to the quantum Dirac operator. This operator acts in the quantum spinor 
bundle, and will be defined in a complete analogy with the classical geometry. We shall study
its properties, including the relations with other important objects of the game. 
In particular, a quantum generalization of the Lichnerowicz formula \cite{li} will be derived. We refer to \cite{ozi2} for a general diagrammatic braided-algebraic foundation of Dirac operator. 

In Section~5 some concluding remarks and observations are made. We shall briefly consider 
the case of arbitrary metric connections (with a possibly non-vanishing torsion), and we shall
also stress specific quantum phenomena appearing in the formalism. As a very instructive 
example, we shall sketch \cite{r-sfera} the construction of the basic objects of the game in the case of the quantum 2-sphere \cite{p} equipped with a canonical spin structure coming from the quantum Hopf fibering. It turns out that in the quantum case the spectrum of the Dirac operator radically differs from the classical situation. The structure group for the Hopf fibering is the classical $\U(1)$ however the differential calculus over it will be {\it quantum}. 

The paper ends with three appendices. 

In the first one, we have included elementary informations about the C*-algebraic
completions of the *-algebras $\cal{V}$ and $\cal{B}$ representing the base space and the bundle. A natural GNS-type
construction is sketched. This construction gives us a natural Hilbert space realization, 
associated to the vertical integration map, of the bundle *-algebra $\cal{B}$. 

In the second appendix, we present in detail the main construction of 
quantum Clifford algebras, associated to general braid operators.
Conceptually, we shall follow \cite{d-ozi}. Our Clifford algebras
will be understood as deformations of the appropriate braided exterior
algebras. As already mentioned, in contrast to \cite{d-ozi} we shall assume
here that the {\it metric components} generate a possibly
non-commutative *-algebra $\mc$ (instead of just being complex numbers). A
special attention will be given to
a discussion of various subtle properties related to the above mentioned
non-commutativity of metric coefficients.

In the last appendix we have collected important properties and definitions related
to the concept of {\it associated vector bundles}. In accordance with \cite{d-tann, d2} these
bundles are defined as $\cal{V}$-bimodules consisting of intertwiners between finite-dimensional
representations of the structure group $G$ and the right action of $G$ on the 
bundle/horizontal forms. It is also explained how to introduce a natural scalar product in the
associated vector bundles, and how our fundamental operators (as Laplacian, Hodge-$*$ and Dirac) naturally act in the intertwiner spaces. 

Finally, a technical remark concerning *-structures. Throughout this paper, *-structures 
on graded algebras will be understood in two ways, depending on the context: as graded-antimultiplicative *-involutions, when dealing with quantum differential forms and exterior 
algebras, or simply as antimultiplicative *-involutions when dealing with braided Clifford algebras and bundles.

\medskip
{\bf Acknowledgements}\par\smallskip
{\it 
I am very grateful to Prof Zbigniew Oziewicz for kindly inviting me to take part in 
the Ixtapa Conference, and for numerous very interesting and fruitful discussions we 
have had on various occasions. }

\section{Quantum Riemannian Geometry}

\subsection{Basic Concepts}

In this section we shall establish the principal geometrical settings for our
consi\-de\-ra\-tions---quantum Riemannian geometry. This will be done 
in two steps. 

At first, we shall recall the definition and basic properties of
quantum principal bundles \cite{d1,d2} and frame structures \cite{d-frm2} on them.  
Quantum frame structures allow us to incorporate into the non-commutative context 
a fundamental concept of coordinate 1-forms, establishing a coherent framework 
for the study of geometrical structures on quantum spaces, in the spirit of classical theory.  
This level is sufficient to introduce general metric connections (together with 
Levi-Civita connections), covariant derivative, curvature and torsion operators. 

Secondly, in order to focus on `true metric spaces', we shall introduce some 
specific analytical properties (that complete our geometrical picture). In
particular, we shall introduce the integration operators for both the frame bundle and the 
base, the Laplace operator, the adjoint differential and the Hodge $*$-operator, 
and prove a couple of important, yet elementary properties. 

Let $G$ be a compact matrix quantum group \cite{w2, w-special}, formally represented by a 
C*-algebra $A$ and a unital *-homomorphism $\phi\colon A\rightarrow A\otimes A$ satisfying
\begin{equation*}
\begin{CD}
A @>{\mbox{$\phi$}}>> A\otimes A\\
@V{\mbox{$\phi$}}VV     @VV{\mbox{$\id\otimes \phi$}}V\\
A\otimes A @>>{\mbox{$\phi\otimes\id$}}> A\otimes A\otimes A
\end{CD}\qquad\qquad
A\otimes A=\lin\overline{\Bigl\{a\phi(b)\Bigr\}}=\lin\overline{\Bigl\{\phi(a)b\Bigr\}}.  
\end{equation*}
The elements of $A$ are interpretable as `continuous functions' over the quantum space 
$G$. The group structure is given by the coproduct map $\phi$. Let $\cal{A}\subseteq A$
be an everywhere dense *-subalgebra corresponding to polynomial functions on $G$. We have 
$\phi(\cal{A})\subseteq\cal{A}\otimes \cal{A}$ and $\cal{A}$ is actually a Hopf *-algebra. 
We shall denote by $\k\colon\cal{A}\rightarrow\cal{A}$ and $\e\colon\cal{A}\rightarrow\Bbb{C}$ the antipode and the counit map respectively. 

Let $M$ be a compact quantum space represented by a *-algebra $\cal{V}$ and let 
$P=(\cal{B}, i, F)$ be a quantum principal $G$-bundle over $M$. By definition \cite{d2},
this means that $\cal{B}$ is a *-algebra,  $i\colon\cal{V}\rightarrow\cal{B}$ is a 
*-monomorphism, and $F\colon\cal{B}\rightarrow \cal{B}\otimes\cal{A}$ is a counital
*-homomorphism such that the following properties hold:

\smallskip
\bla{i}--{\it The action property.} The diagram 
\begin{equation*}\begin{CD}
\cal{B} @>{\mbox{$F$}}>> \cal{B}\otimes\cal{A}\\
@V{\mbox{$F$}}VV     @VV{\mbox{$\id\otimes \phi$}}V\\
\cal{B}\otimes\cal{A} @>>{\mbox{$F\otimes\id$}}> \cal{B}\otimes\cal{A}\otimes\cal{A}
\end{CD}\end{equation*}
is commutative. 

\smallskip
\bla{ii}--{\it The `orbit space' identification.} We have 
$$ i(\cal{V})=\Bigl\{b\in \cal{B}\Bigm\vert F(b)=b\otimes 1\Bigr\}.$$

\smallskip
\bla{iii}--{\it The freeness condition.}  A linear map 
$X\colon\cal{B}\otimes_{\cal{V}}\cal{B}\rightarrow\cal{B}\otimes\cal{A}$ given by 
$$ X(q\otimes b)=qF(b) $$
is surjective. 

\smallskip
It is important to mention that if the map $X$ is surjective then it will be automatically injective, 
so that $X$ is actually bijective \cite{d-ext}. Thus, the freeness condition introduces us 
naturally into the algebraic framework of Hopf-Galois extensions \cite{s-hge}. 

As we shall see later, in a special case of the frame structures representing `vanilla'
quantum Riemannian manifolds, the above freeness condition will be satisfied {\it automatically.} 
 
Now we are going to introduce, following \cite{d-frm2}, the concept of a 
quantum frame structure. 

Let $\Psi$ be a bicovariant \cite{w3} bimodule over $G$. The corresponding
left and right co/action maps will be denoted by $\ell_\Psi\colon\Psi\rightarrow
\cal{A}\otimes\Psi$ and $\rig_\Psi\colon\Psi\rightarrow\Psi\otimes\cal{A}$ 
respectively. Let $\V=\Psi_{\inv}$ be the corresponding left-invariant part.
There exists a natural identification $\Psi\leftrightarrow\cal{A}\otimes\V$ of left 
$\cal{A}$-modules. The structure of $\Psi$ is encoded in the restricted right
action $\v =(\rig_\Psi{\restr}\V)\colon \V\rightarrow \V\otimes\cal{A}$
and a natural right $\cal{A}$-module structure $\circ$ on $\V$, given by
$\vartheta\circ a=\k(a^{(1)})\vartheta a^{(2)}$. 
If $\Psi$ is *-covariant then the space $\V$ is *-invariant. 
We have the following interesting compatibility conditions between $*, \circ$ and $\v$:
\begin{gather*}
\v *=(*\otimes *)\v \qquad (\theta{\circ} a)^*=\theta^*{\circ}\k(a)^*\\
\v (\theta{\circ} a)=\sum_k(\theta_k{\circ} a^{(2)})\otimes\k(a^{(1)})c_k a^{(3)},
\end{gather*}
where $\Sum_k\theta_k\otimes c_k=\v (\theta)$. 

We shall assume that an auxiliary $\v$-invariant scalar product $(\vert)$ is
defined on $\V$. However, for the purposes of our main considerations,
the central role will be played by a {\it noncommutative scalar product} in $\V$,
taking its values in an appropriate *-algebra $\mc$, generated by `abstract
metric tensor coefficients'. 

Let $\tau\colon \V\otimes \V\rightarrow \V\otimes \V$
be the canonical braid operator \cite{w3} associated to $\Psi$. It is
computed in terms of $\v $ and $\circ$, as
$$\tau(\eta\otimes\vartheta)=\sum_k\vartheta_k\otimes(\eta\circ c_k).$$

Let $\V^\wedge$ be the corresponding $\tau$-exterior algebra, obtained from
$\V^\otimes$ by factorizing through the space of quadratic relations 
$\im(I+\tau)$. At this point it is natural to assume that 
$\ker(I+\tau)\neq\{0\}$. This ensures the non-triviality of the higher-order part 
of $\V^\wedge$. 

In what follows, the algebras $\V^\otimes$ and $\V^\wedge$ will be equipped with the 
induced $\circ, *$ and $\v$-structures (these induced structures will be denoted by the
same symbols). The extended structures are constructed by postulating
\begin{gather*}
\v(\vartheta\eta)=\v(\vartheta)\v(\eta)\qquad \v(1)=1\otimes 1\\ 
(\vartheta\eta)^*=(-)^{\partial\vartheta\partial\eta}\eta^*\vartheta^* \\
(\vartheta\eta){\circ} a=(\vartheta{\circ} a^{(1)})(\eta{\circ} a^{(2)})\qquad 1{\circ} a=\e(a)1.
\end{gather*}

The following identities express mutual compatibility between $\tau$ and the maps
$*, \circ$ and $\v$: 
\begin{gather*}
\begin{CD}
\V\otimes\V @>{\mbox{$\v$}}>> \V\otimes\V\otimes\cal{A}\\
@V{\mbox{$\tau$}}VV  @VV{\mbox{$\tau\otimes\id$}}V\\
\V\otimes\V @>>{\mbox{$\v$}}> \V\otimes\V\otimes\cal{A}
\end{CD}\\ 
\tau*=*\tau^{-1}\\
\tau(\psi\circ a)=\tau(\psi)\circ a. 
\end{gather*}

Let us denote by $C_\v\colon\V\rightarrow\V$ the canonical intertwiner
between $\v$ and its second contragradient $\v^{cc}$. By construction, this
map is positive and satisfies 
$$ \v C_\v=(C_\v\otimes \k^2)\v\qquad \tr(C_\v)=\tr(C_\v^{-1}). $$
We shall assume that the scalar product on $\V$ is such that 
\begin{equation}
(x^*,y^*)=(y,C_\v x)\qquad\quad \forall x,y\in\V. 
\end{equation} 
To put it another way around, we can {\it define} $C_\v$ by the above formula. Let us observe 
that this formula implies 
\begin{equation}
*C_\v *=C_\v^{-1} \qquad C_\v=[*]^\dagger *. 
\end{equation}
The operator $C_\v$ is associated to the modular properties of the Haar measure \cite{w2}. 
Polary decomposing the map $*\colon\V\rightarrow\V$ we obtain
\begin{equation}
*=J_\v C_\v^{1/2}=C_\v^{-1/2} J_\v, 
\end{equation}
where $J_\v\colon\V\rightarrow\V$ is an antiunitary involution (in other words 
$J_\v=J_\v^\dagger=J_\v^{-1}$). 

When dealing with various `coordinate expressions', we shall use a fixed  
basis $\{\theta_1,\dots,\theta_d\}$ in $\V$. We shall assume that these vectors satisfy 
$$
(\theta_i,\theta_j)=\delta_{ij}\qquad\quad J_\v(\theta_i)=\theta_i. 
$$
In this basis the representation $\v\colon\V\rightarrow\V\otimes\cal{A}$ is given by 
a unitary matrix $[\v_{ij}]$, so that 
$$
\v(\theta_i)=\sum_j\theta_j\otimes\v_{ji}\qquad\quad C_\v^{1/2}[\v]C_\v^{-1/2}=[\overline{\v}]. 
$$
It is easy to see that the matrix $[C_\v^{1/2}]_{ij}=(\theta_i,C_\v^{1/2}\theta_j)$ is orthogonal. 

As mentioned in the introduction, a `quantum Euclidean' space structure
on $\V$ will be specified by an appropriate
quadratic form $g\colon\V\otimes\V\rightarrow \mc$, playing the role of
the metric, where $\mc$ is a *-algebra generated by matrix elements 
of $g$ and $g^{-1}$, together with a new regular braid operator
$\sigma\colon\V\otimes\V\rightarrow\V\otimes\V$ expressing the twisting
properties of $\mc$ and $\V$. This reflects a fundamental property
of our theory, the non-commutativity of the metric tensor
coefficients and the braided nature of $\mc$. In general $\sigma\neq\tau$. 

The full set of properties involving $g,\sigma$ and $\mc$ is discussed in
Appendix B. Accordingly, we shall also assume that $\{\sigma, \tau\}$
form a pair of coupled braid operators, so that the following natural
identifications
\begin{equation}\label{antisym-id}
\Bigl\{\text{$\tau$-antisymmetric $n$-tensors}\Bigr\}\leftrightarrow
\im(A_\sigma^n)
\end{equation}
hold for each $n\geq 2$.
Here $A_\sigma^n\colon\V^{\otimes n}\rightarrow\V^{\otimes n}$ are the braided
$\sigma$-antisymmetrizers. By definition, $\tau$-antisymmetric tensors
change the sign under all the corresponding $\tau$-transpositions.
Our $\tau$-exterior algebra will be realizable as follows 
$$\V^\wedge\leftrightarrow
\V^\otimes/\ker(A_\sigma)\leftrightarrow \im(A_\sigma).$$
These identifications are constructed with the help of $A_\sigma$. 
It is worth mentioning that the grading and the *-structure are preserved
in this picture. The maps $C_\v$ and $J_\v$ will be extended
(by multiplicativity/anti) to $\V^\otimes$. The space $\V^\wedge$ is
invariant under the action of these maps. As discussed in Appendix B, there is a 
natural twisted tensor product of algebras  $\V^\wedge$ and $\mc$, and in 
such a way we obtain an extended braided exterior algebra $\Vmc^\wedge$. This
is a braided exterior algebra built over $\Vmc$ with the help of the extended
braid $\sigma\colon\Vmc\otimes_\mc\Vmc\rightarrow \Vmc\otimes_\mc\Vmc$. 

\medskip
********\par
By definition, a {\it frame structure} on a quantum principal bundle $P$ is given by a special 
graded *-algebra $\hor_P$, equipped with a first-order hermitian antiderivation
$\nabla\colon\hor_P\rightarrow\hor_P$. The algebra $\hor_P$ is defined as 
$\hor_P\leftrightarrow\cal{B}\otimes\V^\wedge$ at the level of vector spaces, while the
product and the *-structure are given by formulae
\begin{align*}
(q\otimes\vartheta)(b\otimes\eta)&=\sum_kqb_k\otimes(\vartheta\circ c_k)
\eta\\
(b\otimes\vartheta)^*&=\sum_k b_k^*\otimes(\vartheta^*\circ c_k^*)
\end{align*}
where $\Sum_k b_k\otimes c_k=F(b)$. The elements of $\hor_P$ are interpreted as 
`quantum horizontal forms'. The elements of $\V^\wedge$, viewed in the framework of $\hor_P$, 
are interpretable as analogs of natural `coordinate forms' in classical theory of frame bundles. 
We see that $\hor_P^0=\cal{B}$. The maps $F$ and $\v\colon\V^\wedge\rightarrow\V^\wedge
\otimes\cal{A}$ naturally combine to a unital *-homomorphism $\Fh\colon\hor_P\rightarrow\hor_P
\otimes\cal{A}$ satisfying
\begin{gather*}
(\id\otimes \phi)\Fh=(\Fh\otimes\id)\Fh\\
(\id\otimes\e)\Fh=\id.
 \end{gather*}
The map $\Fh$ plays the role of the right action of $G$ on horizontal forms. The corresponding
$\Fh$-fixed-point graded *-subalgebra $\Omega_M\subseteq\hor_P$ plays the role of the 
differential forms on the base manifold $M$. Accordingly $\Omega_M^0=\cal{V}$. 

The map $\nabla\colon\hor_P\rightarrow\hor_P$ completes the picture of a frame bundle. It 
corresponds to the Levi-Civita connection in classical theory. By definition, $\nabla$ intertwines
the action $\Fh$ and satisfies $$\nabla\bigl\{\nabla(\cal{V})\V\bigr\}=\{0\}.$$ In particular 
$\nabla$ vanishes on the subalgebra generated by $\nabla(\cal{V})$ and $\V^\wedge$. 

Finally, it is assumed that there exist linear maps $b_{\alpha}:\theta\mapsto 
b_{\alpha}(\theta)\in\cal{B}$ and elements $f_\alpha\in
\cal{V}$ satisfying
\begin{gather}
1\otimes\theta=\sum_\alpha b_{\alpha}(\theta)\nabla(f_\alpha)\label{frm-rel1}\\
F[b_{\alpha}(\theta)]=(b_{\alpha}\otimes \id)\v(\theta).\label{frm-rel2} 
\end{gather}
This is a kind of a completeness condition. From the above mentioned postulates, 
it follows that $\nabla$ is reduced in $\Omega_M$, and that the restriction map
$\dM\colon\Omega_M\rightarrow\Omega_M$ is a hermitian differential (corresponding
to the standard exterior derivative of differential forms). Moreover, it can be shown that
the algebra $\Omega_M$ is generated by $\cal{V}$ and $\dM(\cal{V})$. 

We can introduce an explicit `coordinate' description of $\nabla$, by the formula
\begin{equation}\label{nabla-coor}
\nabla(b)=\sum_i \partial_i(b)\otimes \theta_i,
\end{equation}
where $b\in\cal{B}$. The maps $\partial_i\colon\cal{B}\rightarrow\cal{B}$ are counterparts of 
canonical horizontal coordinate vectors fields. They completely determine $\nabla$. 
In particular, the $\Fh$-covariance of $\nabla$ is equivalent to the property
\begin{equation}\label{frm-rel3}
F\partial_i(b)=\sum_{jk} \partial_j(b_k)\otimes c_k\k^{-1}(\v_{ij}),
\end{equation}
where $F(b)=\Sum_k b_k\otimes c_k=F(b)$ and $\Sum_j\theta_j\otimes 
\v_{ji}=\v(\theta_i)$. The graded Leibniz rule for $\nabla$ is translated into the 
system of equations
\begin{equation}\label{X-leibniz}
\partial_i(qb)=q\partial_i(b)+\sum_j\partial_j(q)\mu_{ji}(b),
\end{equation}
where $\mu\colon\cal{B}\rightarrow \mathrm{M}_d(\cal{B})$ is a monomorphism defining
the right $\cal{B}$-module structure on $\hor_P$. In other words, 
$$
\theta_i b=\sum_j \mu_{ij}(b)\theta_j. 
$$
The maps $\mu_{ij}\colon\cal{B}\rightarrow\cal{B}$ are expressible via the 
right $\cal{A}$-module structure on $\V$ and the map $F$ as 
\begin{equation}\label{mu-F-nu}
\mu_{ij}(b)=\sum_k b_k \nu_{ij}(c_k),
\end{equation}
where $\nu\colon\cal{A}\rightarrow \mathrm{M}_d(\Bbb{C})$ is a unital homomorphism given
by 
\begin{equation}
\theta_i{\circ} a=\sum_j \nu_{ij}(a)\theta_j, \qquad\quad \mu(f)=fI_d \quad f\in\cal{V}.  
\end{equation}
The introduced coordinate vector fields fit into a general framework for quantum vector
fields, introduced in \cite{bor}. 

The following identity expresses the compatibility between $\circ$ and the *-structure 
on $\V$: 
\begin{equation}\label{nu-Cv}
\nu[\k(a)^*]=C^{-1/2}_\v\overline{\nu(a)}C^{1/2}_\v. 
\end{equation}
The compatibility between $\mu$ and the *-structure on $\hor_P$ reads
\begin{equation}\label{*-mu}
\delta_{ij}b^*=\sum_k\bigl\{C^{1/2}_\v\mu C^{-1/2}_\v\bigr\}_{kj}\bigl\{\mu_{ik}(b)^*\bigr\}. 
\end{equation}

Along these lines, it is also possible to express the hermiticity property of $\nabla$ in terms of 
the maps $\partial_i$---a direct computation gives
\begin{equation}\label{X*}
\partial_i=\sum_j[C^{1/2}_\v\mu]_{ji}\{\partial^*_j\} \qquad \quad\partial_j^*=*\partial_i{*}. 
\end{equation}
The structure group $G$ corresponds to a transformation group of `local orthonormal 
frames'.  The frame structure allows us to think of $P$ as the bundle of `orthonormal frames' 
over $M$. In accordance with this analogy, it would be natural to assume that $\v$ is faithful (in
other words $G$ is completely determined by its action on $\V$, which means that $\cal{A}$ is
generated by the matrix elements of $\v$). However, from the point of view of our spinorial 
constructions it is natural to allow the situations where $\v$ is not faithful
($G \leftrightarrow \mathrm{Spin}(m)$). This allows us
to include the spin structures within the framework of the frame structures. 

\begin{lem}
In the case of `truly frame bundles' when $\v$ is faithful the freeness condition for $F$ 
is automatically fulfilled. 
\end{lem}

\begin{pf} Let us assume that $\v$ is faithful and that all the conditions on the frame structure 
are satisfied, except possibly the freeness of $F$. Let us first observe that 
$$
\sum_\alpha b_{\alpha j}F(q_{\alpha i})=1\otimes \k^{-1}(\v_{ij})
$$
as it follows from \eqref{frm-rel1} and \eqref{frm-rel3}. Here
$q_{\alpha i}=\partial_i(f_\alpha)$ and $b_{\alpha i}=b_{\alpha}(\theta_i)$. 

This implies that all the elements of the form $b\otimes \v_{ij}$ are in the image of the 
canonical map $X\colon\cal{B}\otimes_{\cal{V}}\cal{B}\rightarrow\cal{B}\otimes\cal{A}$. 
Now playing with the elementary \cite{d2} properties of $X$ and the fact that $\cal{A}$ is 
generated by $\v_{ij}$, we conclude that $X$ is surjective. 
\end{pf}

We are now going to write down very important commutation relations between
coordinate vector fields $\partial_i$, involving the curvature tensor. Let us recall that
the curvature $\varrho_\nabla\colon\cal{A}\rightarrow \zh_P^2$ of $\nabla$ is uniquely 
determined through a fundamental identity
\begin{equation}\label{nabla2}
\nabla^2(b)=-\sum_k b_k\varrho_\nabla(c_k). 
\end{equation}
Here $\zh_P$ is the graded commutant of $\Omega_M$ in $\hor_P$. The fact that the
curvature always take values from $\zh_P^2$ implies strong constraints for possible
forms of a Levi-Civita connection in non-commutative geometry. 

The canonical inclusion of the exterior algebra into the tensor algebra allows us to define the
`components' $\varrho_\nabla^{ij}\colon\cal{A}\rightarrow\cal{B}$ of the curvature, by the
formula 
\begin{equation}
\varrho_\nabla(a)=\frac{1}{2}\sum_{ij}\varrho_\nabla^{ij}(a)\otimes \{\theta_i\otimes\theta_j\}. 
\end{equation}
The components by construction satisfy a $\tau$-antisymmetricity relations
\begin{equation}
\varrho_\nabla^{ij}=-\sum_{kl}\tau_{kl}^{ij}\varrho_\nabla^{kl}
\end{equation}
where $\tau(\theta_i\otimes\theta_j)=\Sum_{kl}\tau^{kl}_{ij}\theta_k\otimes\theta_l$. 

\begin{lem} The following identity holds:
\begin{equation}\label{com-rel-X}
\partial_i\partial_j(b)-\sum_{kl}\sigma^{ij}_{kl}\partial_k\partial_l(b)+\frac{1}{2}\sum_\alpha b_\alpha\varrho_\nabla^{ij}(c_\alpha)=0, 
\end{equation}
where $\Sum_{kl}\sigma_{ij}^{kl}\theta_k\otimes\theta_l=\sigma(\theta_i\otimes\theta_j)$.  
\end{lem}

\begin{pf}
This identity follows by rewriting \eqref{nabla2} in the coordinate form,
using the definitions of $\tau$ and $\sigma$, and the above coordinate expression for the
curvature.
\end{pf}

\subsection{Integration Operators}

Let us start by formulating some additional geometrical assumptions about the 
quantum base space $M$. At first, we shall assume that $\cal{V}$
is realized as an everywhere dense *-subalgebra of a unital C*-algebra $\cstV$. We shall
also assume that a faithful state $\omega_M\colon\cstV\rightarrow \Bbb{C}$ is given
(representing a `measure' on $M$). By definition, the faithfulness property means that
$\omega_M$ is {\it strictly positive} on the positive elements. Furthermore, we shall assume that
$\omega_M$ admits a modular operator $\modM\colon\cal{V}\rightarrow\cal{V}$
of the form
\begin{equation}
\omega_M(fg)=\omega_M(\modM(g)f) \qquad\quad\forall f,g\in\cal{V}. 
\end{equation}
The modular operator $\modM$ is uniquely determined by the state $\omega_M$. The following
identities hold
\begin{gather*}
\modM(fg)=\modM(f)\modM(g)\\
\modM {*}\modM={*}
\end{gather*}
and in particular we see that $\modM$ is necessarily bijective. 

Next, we shall introduce compatibility conditions between the quantum principal bundle 
$P$ with the corresponding frame structure, and the measure $\omega_M$ with
the associated modular automorphism $\modM$. 
 
Let us introduce a `vertical integration' operator $\vint\colon\cal{B}\rightarrow\cal{V}$ by
\begin{equation}
i[\vint(b)]=(\id\otimes h)F(b),
\end{equation}
where $h\colon\cal{A}\rightarrow \Bbb{C}$ is the Haar measure \cite{w2} of $G$. 
Our first assumption on $P$ will be:

\smallskip {\samepage 
{\bf \begin{flushright} The Strict Positivity\end{flushright}}

\smallskip\noindent
{\it The map $\vint$ is strictly positive. In other words, 
$$ \vint(b^*b)\geq 0 \qquad\quad\forall b\in\cal{B} $$
and if $\vint(b^*b)=0$ then $b=0$. } }

\smallskip
This property ensures that the *-algebra $\cal{B}$ is closable into a C*-algebra $\cstB$, 
using a natural GNS-type faithful *-representation by bounded operators.
The representation is constructed with the help of $\vint$ and $\omega_M$, as
discussed in Appendix A. By combining $\vint$ and $\omega_M$ we arrive 
to a faithful state $\omega_P\colon\cstB\rightarrow\Bbb{C}$ playing the role of
the `measure' on the bundle $P$. Let us assume that $\omega_P$ admits a modular operator 
$\modP\colon\cal{B}\rightarrow\cal{B}$. 

\begin{lem} We have 
\begin{equation}\label{FmodP}
F\modP=(\modP\otimes \k^2)F. 
\end{equation}
In particular $\modP[\cal{V}]=\cal{V}$ and $\modP\restr\cal{V}=\modM$. 
\end{lem}

\begin{pf} The invariance property of $\omega_P$ can be written in the form 
$$ \sum_l\omega_P(bq_l)\otimes\k(d_l)=\sum_k\omega_P(b_kq)\otimes c_k. $$ 
This, together with the definition of $\modP$ gives
\begin{multline*}
\sum_l \omega_P\bigl[\modP(q_l)b\bigr]\otimes\k(d_l)=\sum_l \omega_P(bq_l)\otimes\k(d_l)=
\sum_k\omega_P(b_kq)\otimes c_k\\=\sum_k\omega_P(\modP(q)b_k)\otimes c_k=
\omega_P(z_jb)\otimes \k^{-1}(e_j), 
\end{multline*}
where $F\modP(q)=\Sum_j z_j\otimes e_j$. Hence, it follows
that $F\modP(q)=\Sum_l\modP(q_l)\otimes \k^2(d_l)$. 
\end{pf}
 
The space $\cal{B}$ is equipped with a natural $F$-invariant scalar product given by 
$$ <b,q>=\omega_P(b^*q). $$
Playing with the definitions of $\modP$ and $<,>$ it follows that 
$*\colon\cal{B}\rightarrow\cal{B}$ is formally adjointable. Accordingly, 
\begin{gather*}
<b,q^*>=<q,[*]^\dagger(b)>\qquad\quad\forall q,b\in\cal{B}\\
[*]^\dagger=\modP^{-1}*=*\modP. 
\end{gather*}

Let $\cal{T}$ be a complete set of mutually non-equivalent irreducible representations 
of $G$. By decomposing $\cal{B}$ into the multiple irreducible submodules relative to
the action $F$, we arrive to 
$$\cal{B}=\sideset{}{^\oplus}\sum_{\alpha\in \cal{T}} \cal{B}^\alpha 
\qquad\cal{B}^\alpha\leftrightarrow\bim{\alpha}\otimes H_\alpha.$$
Here we have used intertwiner bimodules $\bim{\alpha}=\Mor(\alpha, F)$, and $H_\alpha$
are the corresponding representation spaces. If $\varnothing$ is the trivial 
representation in $\Bbb{C}$ then 
$\cal{B}^\varnothing\leftrightarrow\cal{V}$. For each $\alpha\in\cal{T}$ let us denote by
$\{\,\}_\alpha\colon\cal{B}\rightarrow\cal{B}^\alpha$ the corresponding projection map. 
The vertical integration map is given by 
projecting $\cal{B}$ onto $\cal{V}$ in other words $\vint\leftrightarrow\{\,\}_\varnothing$. 
Our second extra assumption on $P$ will be:

\smallskip
{\bf \begin{flushright} The Horizontal Homogeneity\end{flushright}}

\smallskip\noindent
{\it Let  us consider the elements $b_{\alpha i}=b_\alpha(\theta_i)$. Then we have} 
\begin{equation}                   
\omega_M\Bigl\{\sum_{ij} [C_{\v}^{-1}]_{ji}\partial_i(b_{\alpha j}f)\Bigr\}=0, \qquad\quad 
\forall f\in\cal{V}. 
\end{equation}

This condition expresses the idea that the measure on $M$ is `homogeneous'. As we have
seen, the maps $\partial_i\colon\cal{B}\rightarrow\cal{B}$ play the role of canonical horizontal 
vector fields, and $\omega_M$ should be invariant under the appropriate 
`infinitesimal horizontal transformations'. 
Let us observe that first-order `differential' operators 
$T_\alpha\colon\cal{V}\rightarrow\cal{V}$ defined by 
$$T_\alpha(f)=\sum_{ij} [C_{\v}^{-1}]_{ji}\partial_i(b_{\alpha j}f) $$
figure in the above expression.
These operators are naturally associated to the frame structure. 

\begin{lem} Under the above homogeneity and positivity assumptions, we have
\begin{equation}\label{wP-X}
\omega_P\bigl[\partial_i(b)\bigr]=0 \qquad\forall b\in\cal{B}. 
\end{equation} 
\end{lem}

\begin{pf} Let us observe that \eqref{wP-X} is non-trivial only for the 
$F$-invariant component $\{\partial_i(b)\}_\varnothing$---all other components of 
$\partial_i(b)$ are already annihilated by the vertical integration map $\vint$. On the 
other hand, it is easy to see that 
$$
\{\partial_i(b)\}_\varnothing=\Sum_j\partial_j(q_j), 
$$
where the elements $q_j$ form an appropriate $F$-multiplet
$$
F(q_j)=\sum_k q_k\otimes \k^{-2}(u_{kj}). 
$$
All such $F$-multiplets are of the form 
$$
q_j=\sum_{i\alpha} [C^{-1}_\v]_{ij}b_{\alpha i} f_\alpha
$$
where $f_\alpha\in\cal{V}$ are arbitrary. Now the statement of the lemma follows from the 
homogeneity assumption. 
\end{pf}

It is interesting to calculate the formally adjoint operators for 
important coordinate maps $\mu_{ij}\colon\cal{B}\rightarrow\cal{B}$ and 
$\partial_i\colon\cal{B}\rightarrow\cal{B}$. 

\begin{lem} The maps $\mu_{ij}$ and $\partial_i$ are formally adjointable and
\begin{gather}
\mu^\dagger(b)=C^{1/2}_\v\mu(b)C^{-1/2}_\v\label{mu-adj}\\
-\partial_i^\dagger=\sum_j[C^{-1/2}_\v]_{ji}\partial_j=\sum_j
\mu^\dagger_{ji}\{\partial^*_j\}.\label{X-adj} 
\end{gather} 
\end{lem}

\begin{pf} Let us first calculate the adjoint of $\mu_{ij}$. Applying \eqref{mu-F-nu} and 
\eqref{nu-Cv}, and
the $F$-invariance of the scalar product $<,>$ we find
\begin{multline*}
<q,\mu_{ij}(b)>=\sum_k<q,b_k>\nu_{ij}(c_k)=\sum_l<q_l,b>\nu_{ij}[\k(d_l)^*]\\
=\sum_l<q_l,b>\Bigl\{C^{-1/2}_\v\overline{\nu(d_l)}C^{1/2}_\v \Bigr\}_{ij}=
\sum_l<q_l\Bigl\{C^{1/2}_\v\nu(d_l)C^{-1/2}_\v \Bigr\}_{ij},b>\\
=<\Bigl\{C^{1/2}_\v\mu(q)C^{-1/2}_\v \Bigr\}_{ij},b>
\end{multline*}
which proves \eqref{mu-adj}. Next, applying \eqref{wP-X} and \eqref{X-leibniz} we obtain 
\begin{equation*}
\begin{split}
0=\omega_P \partial_i(q^* b)&=\omega_P\Bigl\{q^*\partial_i(b)\Bigr\}+\sum_j\omega_P
\Bigl\{\partial_j(q^*)\mu_{ji}(b)\Bigr\}=\\
{}&=<q,\partial_i(b)>+\sum_j<\partial^*_j(q),\mu_{ji}(b)> 
\end{split}
\end{equation*}
and thus the second equality in \eqref{X-adj} holds. The first equality in \eqref{X-adj} directly 
follows from the second, together with \eqref{mu-adj} and \eqref{X*}. 
\end{pf}

We are going to construct the integration operator for horizontal forms. This will be done
by combining the measure $\omega_P$ and a `coordinate volume form' the existence of 
which will be ensured by our next extra condition:

\smallskip
{\bf \begin{flushright}Self-Duality of Coordinate Forms\end{flushright}}

\smallskip\noindent
{\it There exists a number $m\in\Bbb{N}$ such that $\V^{\wedge m}\leftrightarrow\Bbb{C}$
and $\V^{\wedge k}=\{0\}$ for each $k>m$. }

\smallskip
In other words, we can introduce the `volume element' as a single generator 
$w=w^*\in\V^{\wedge m}$. 

Actually, the above condition follows from a simple assumption that the braided 
exterior algebra $\V^\wedge$ is finite-dimensional. To see this, let us consider 
an auxiliary braided Clifford algebra $\cal{D}$ associated to $\{\V,\V^*,\sigma\}$. Here
it is assumed that $\sigma$ acts on all possible tensor products involving 
$\V$ and $\V^*$. Braided exterior algebras $\V^\wedge$ and $\V^{\wedge*}=[\V^*]^\wedge$
are subalgebras of $\cal{D}$ and we have the following natural vector space identifications
$$ \cal{D}\leftrightarrow \V^\wedge\otimes \V^{\wedge *}\leftrightarrow \V^{\wedge *}\otimes \V^\wedge$$
induced by the product map. The algebra $\cal{D}$ naturally acts on the space 
$\V^\wedge$---the action is constructed with the help of the canonical duality between $\V^\wedge$ and $\V^{\wedge*}=[\V^*]^\wedge$. Explicitly, we have 
$$ \V^\wedge=\cal{D}\otimes_{*\wedge}\Bbb{C}$$
at the level of left $\cal{D}$-modules. In a similar manner, we can construct the dual 
module, by considering
$$ \V^{\wedge*}=\cal{D}\otimes_\wedge\Bbb{C}$$ 
and in both formulas $\Bbb{C}$ is viewed as a trivial module over the corresponding
exterior algebra. Both modules are irreducible. Now if $\V^\wedge$ is finite-dimensional, then $\cal{D}$ will be the full endomorphism algebra, and as such it will be allowed 
to have only one irreducible representation, up to equivalence. 
In particular, there exists a bijective 
$\cal{D}$-linear map $Y\colon\V^\wedge\rightarrow \V^{\wedge *}$.  

This map is a blueprint for the Hodge *-operator. It intertwines contractions and wedge-products, and in particular it follows that $Y[\V^{\wedge k}]=\V^{\wedge m-k*}$ where 
$m\in\Bbb{N}$ is the dimension giving the volume element. The formula
\begin{equation}
j(x,y)w=x\wedge y
\end{equation}
defines a nondegenerate pairing $j\colon\V^{\wedge k}\times\V^{\wedge m-k}\rightarrow\Bbb{C}$. 
Moreover, it follows that there exists a unique grade-preserving map $\ll\colon\V^\wedge\rightarrow
\V^\wedge$ such that 
\begin{equation}
j(y,x)=(-)^{\partial x\partial y}j(\ll(x),y)\qquad\qquad\forall x,y\in\V^\wedge. 
\end{equation}

From the definition of $w$ it follows that 
\begin{equation}\label{Q-l}
\v(w)=w\otimes Q\qquad\qquad w{\circ}a=\lambda(a) w, 
\end{equation}
where $Q\in\cal{A}$ is a `quantum determinant' such that
$$ \phi(Q)=Q\otimes Q, \quad\qquad \k(Q)=Q^{-1}=Q=Q^* $$
and $\lambda\colon\cal{A}\rightarrow \Bbb{C}$ is a unital multiplicative functional
satisfying 
$$ \overline{\lambda(a)}=\lambda[\k(a)^*].$$
It is easy to see that the introduced objects satisfy
\begin{gather*}
\ll(x\circ a)=\ll(x){\circ}\k^2(a)\\
\v\ll=(\ll\otimes\id)\v\\
Qa=\ad(a,\lambda)Q\\ 
j(x{\circ} a^{(1)},y{\circ}a^{(2)})=\lambda(a)j(x,y)\\
*\ll*=\ll^{-1}\qquad \overline{j(x,y)}=(-)^{\partial x\partial y}j(y^*,x^*).
\end{gather*} 
Let us observe that in general $m\neq d$, in contrast with classical
geometry.

Now we shall introduce the integration map $\intP\colon\hor_P\rightarrow \Bbb{C}$. 
The definition is straightforward
\begin{equation}\label{intP-def}
\intP(b\otimes \vartheta)=\begin{cases}\omega_P(b)&\text{for $\vartheta=w$,}\\
0 & \text{if $\deg \vartheta< m$}.\end{cases}
\end{equation}

\begin{lem}  We have
\begin{equation}\label{intP-nabla}
\intP\nabla(\varphi)=0\qquad\qquad\forall\varphi\in \hor_P. 
\end{equation}
Moreover, 
\begin{equation} 
\intP[\varphi^*]=(\intP\varphi)^*\qquad\sum_k(\intP\varphi_k)\otimes c_k=(\intP\varphi)\otimes Q,
\end{equation}
where $\Sum_k\varphi_k\otimes c_k=\Fh(\varphi)$. 
\end{lem}

\begin{pf} Property \eqref{intP-nabla} follows directly from \eqref{wP-X} and the definition 
of $\intP$. Hermiticity and covariance of $\intP$ are simple consequences of the similar properties for
the measure $\omega_P$. 
\end{pf}

Let us now observe that $Q$, as a hermitian involution, naturally decomposes into hermitian
projections $Q_+,Q_-$ where $Q_\pm=1/2\pm Q/2$. If the map 
$\lambda\colon\cal{A}\rightarrow\Bbb{C}$ is in addition {\it central}, then it will be possible
to pass to the corresponding `components' of $G$, determined by $Q_\pm$. In particular, we can 
factorize through the Hopf *-ideal generated by $Q_-$, reducing to the case $Q=1$. This is the
quantum version of {\it unimodularity}---passing from $\O(d)$ to $\SO(d)$ groups. 

Geometrically, such a restriction means that we are dealing with {\it oriented manifolds}. 
In what follows it will be assumed that orientability property holds (and centrality of $\lambda$, as 
a necessary consistency condition). 

For the end of this subsection, we shall analyze modular properties of the integration map. 
\begin{lem} The following identity holds
\begin{equation}\label{intP-mod}
\intP \!\varphi\psi=(-)^{\partial\varphi\partial\psi}\intP \!\modH(\psi)\varphi
\end{equation}
where $\modH\colon\hor_P\rightarrow\hor_P$ is a grade-preserving homomorphism defined 
by $\modH(b\otimes\vartheta)=\modP(b_\lambda)\otimes\ll(\vartheta)$ and 
$[]_\lambda=(\id\otimes\lambda)F$. 
\end{lem}

\begin{pf} Assuming that $\varphi=b\otimes \vartheta$ and $\psi=q\otimes\eta$ and performing
elementary transformations we find
\begin{multline*}
\intP \!\varphi\psi=\sum_l\intP\Bigl\{ bq_l\otimes(\vartheta{\circ}d_l)\eta\Bigr\}=
\sum_l\omega_P(bq_l)j(\vartheta{\circ}d_l,\eta)\\=\sum_k\omega_P(\modP(q)b_k)
j(\vartheta{\circ}\k^{-1}(c_k),\eta)=\sum_k(-)^{\partial\varphi\partial\psi}\omega_P(\modP(q)b_k)
j(\ll(\eta),\vartheta{\circ}\k^{-1}(c_k))\\=
\sum_k(-)^{\partial\varphi\partial\psi}\omega_P(\modP(q)b_k)j(\ll(\eta){\circ}c_k^{(1)}\lambda^{-1}(c_k^{(2)}),\vartheta)\\
=\sum_k(-)^{\partial\varphi\partial\psi}\omega_P
(\modP(q_\lambda)b_k)j(\ll(\eta){\circ}c_k,\vartheta)
=(-)^{\partial\varphi\partial\psi}\intP\!\Bigl\{\modH(\psi)\varphi\Bigr\}. 
\end{multline*}
Let us also observe that 
$$
\Fh\modH=(\modH\otimes\k^2)\Fh\qquad {*}\modH{*}=\modH^{-1}.  
$$
The map $\modH$ is uniquely determined by \eqref{intP-mod}. 
\end{pf}

Through this paper we shall make an extensive use of an {\it extended horizontal forms} 
algebra $\horPmc$ obtained by mixing $\mc$ with the standard horizontal forms. More precisely, 
$\horPmc$ is obtained by taking the twisted tensor product between $\cal{B}$ and the extended braided exterior algebra $\Vmc^\wedge$. 
We have obviously natural left/right $\cal{B},\mc$-module identifications 
$$ \horPmc\leftrightarrow \mc\otimes\hor_P\leftrightarrow \hor_P\otimes \mc. $$
The action $\Fh$ naturally extends, with the help of $\vmc\colon\mc\rightarrow\mc\otimes \cal{A}$, to the action $\Fh\colon\horPmc\rightarrow\horPmc\otimes\cal{A}$. 
The frame structure $\nabla$ naturally extends, by $\mc$-left/right linearity, to a hermitian antiderivation $\nabla\colon\horPmc\rightarrow\horPmc$. 

The base space algebra $\Omega_M$ is naturally included in $\Wmc$---which is defined 
as the $\Fh$-fixed point subalgebra of $\horPmc$. We shall freely pass from extended to 
non-extended objects, and vice versa. 

\subsection{The Hodge Operator---Linear Algebra}

We shall first introduce the Hodge $*$-operator on $\Vmc^\wedge$. Then it will be extended to 
the level of $\horPmc$. We shall assume that $\Vmc^\wedge$ is
equipped with a $\mc$-valued quadratic form $g_\wedge$ and the associated
scalar product $\langle\rangle$, as explained in Appendix~B. We have
$$ \Vmc^{\wedge m}=\mc w=w \mc \qquad w \alpha =S(\alpha) w$$
where $S\colon\mc\rightarrow \mc$ is an automorphism satisfying
\begin{gather*}
S(\alpha){\circ}a=S(\alpha{\circ}a) \qquad *S*=S^{-1}\\
\vmc S=(S\otimes\id)\vmc. 
\end{gather*}

The map $j$ is straightforwardly extendible to a $\mc$-valued pairing
acting within $\Vmc^\wedge$. We have
\begin{gather*}
j(\varphi,\psi q)=j(\varphi,\psi)S(q)\qquad j(q\varphi,\psi)=
q j(\varphi,\psi)\\
j(\varphi,q \psi)=j(\varphi q,\psi)\\ 
j(\varphi{\circ} a^{(1)}, \psi{\circ} a^{(2)})=\lambda(a^{(1)})j(\varphi,
\psi){\circ}a^{(2)}
\end{gather*}

\begin{pro} \bla{i} The formula
\begin{equation}\label{hstr-def}
g_\wedge(x,y)=j (x, \hstr[y]) 
\end{equation}
uniquely defines a linear operator $\hstr\colon\Vmc^\wedge\rightarrow\Vmc^\wedge$ such that
$\hstr(\Vmc^{\wedge k})\subseteq\Vmc^{\wedge m-k}$. 

\smallskip
\bla{ii} The map $\hstr$ is bijective, and we have 
\begin{gather}
\hstr[x q]=\hstr[x]S^{-1}(q)\qquad \hstr[q x]=q \hstr[x] \qquad q\in \mc\\
\v\hstr=(\hstr\otimes\id)\v\label{v-hstr}\\
\hstr(\vartheta\circ a_\lambda)=\hstr(\vartheta)\circ a.\label{circ-hstr} 
\end{gather}
\end{pro}

\begin{pf} The $\v$-covariance of $\hstr$ follows from 
formula \eqref{hstr-def} and the covariance of all defining entities. 
The $\mc$-compatibility follows from a similar property for $j$.
Finally, performing elementary transformations we obtain
\begin{multline*}
g_\wedge(\xi,\vartheta){\circ a}=g_\wedge(\xi{\circ} a^{(1)},
\vartheta{\circ} a^{(2)})=j(\xi{\circ} a^{(1)},
\hstr[\vartheta{\circ} a^{(2)}])
=j(\xi,\hstr[\vartheta]){\circ} a\\=j(\xi{\circ} a^{(1)},
\hstr[\vartheta]{\circ} a^{(2)})\lambda^{-1}(a^{(3)}),
\end{multline*} 
and hence \eqref{circ-hstr} holds.
\end{pf}

The automorphism $S\colon\mc\rightarrow \mc$ has a simple structure, as it is
sufficient to calculate its action on the elements of the form
$g(x,y)$ where $x,y\in \V$.

By construction, we have $\sigma\colon w\otimes x\mapsto T(x)\otimes w$ and $\sigma \colon
y\otimes w\mapsto w\otimes T^*(y)$ where $T\colon \V\rightarrow \V$ is a bijective
linear operator. Using this, and iteratively applying the definition of
the $\mc$-bimodule structure on $\Vmc$ we find 
\begin{equation}\label{S-<>}
S\bigl\{\langle x ,y\rangle\bigr\}=\langle T^{-1}(x), T(y)\rangle. 
\end{equation}
The map $T$ naturally extends to a unital automorphism $T\colon\V^\wedge
\rightarrow \V^\wedge$, and furthermore to $T\colon\Vmc^\wedge
\rightarrow \Vmc^\wedge$ by imposing $T(xq)=T(x)S(q)$ and $T(qy)=S(q)T(y)$.   
The above formula remains valid for arbitrary elements of our braided exterior algebra. 

\begin{lem} We have
\begin{equation}\label{hstr-adj}
\hstr^\dagger=\hstr, \qquad\qquad T^\dagger=T.  
\end{equation}
In other words, 
$$
\langle x ,T(y)\rangle=S\bigl\{\langle T(x),y\rangle\bigr\} \qquad 
\langle x ,\hstr(y)\rangle=S^{-1}\bigl\{\langle \hstr(x), y\rangle\bigr\}
$$
Here, the concept of the adjoint operator is appropriately $S$-twisted---a necessary
consistency condition, having in mind a right $S$-twisted $\mc$-linearity of $\hstr$ and $T$. 
\end{lem}

\begin{pf} In fact, equation \eqref{S-<>} expresses the selfadjointness of $T$. We compute 
$$
[x^*y]_m^*=w \langle x,\hstr^{-1}(y)\rangle^* =S\bigl\{\langle
\hstr^{-1}(y), x\rangle\bigr\} w=
[y^*x]_m=\langle y, \hstr^{-1}(x)\rangle w. 
$$
This shows that $\hstr$ is selfadjoint. 
\end{pf}

We conclude this subsection by connecting the contraction and multiplication operators in 
$\Vmc^\wedge$. By definition, contraction operators $\iota[x]\colon \Vmc^\wedge
\rightarrow\Vmc^\wedge$ are given by 
\begin{equation} 
\iota[x]\psi=(g\otimes\id^{n-1})(x\otimes\psi)\qquad\forall\psi\in \Vmc^{\wedge n},
\end{equation}
where $\psi$ are realized in the tensor algebra. 

It is easy to see that these operators are $\sigma$-braided antiderivations
(they satisfy the $\sigma$-braided Leibniz rule). In other words,
\begin{equation}\label{contr-br}
\iota[x]y+\sum_\alpha y_\alpha\iota[x_\alpha]=g(x,y) \qquad\forall x,y\in\Vmc
\end{equation}
where $\Sum_\alpha y_\alpha\otimes x_\alpha=\sigma(x\otimes y)$. 

Furthermore, using a natural $\mc$-valued scalar product in $\Vmc^\wedge$ it
follows that
\begin{equation}\label{contr-adj}
[x\wedge()]^\dagger=\iota[x^*] \qquad \forall x\in \Vmc. 
\end{equation}
In other words, the contraction operators are the adjoint maps of the 
corresponding multiplication operators. 

\begin{lem}The following identity holds
\begin{equation}\label{e-adj}
\iota[e]=\hstr^{-1}[e\wedge()]\hstr. 
\end{equation}
\end{lem}

\begin{pf} We have
$$
[x^*\hstr\iota[e](y)]_m=\langle x,\iota[e](y)\rangle w=
\langle e^*\wedge x,y\rangle w=[x^*\{e\wedge\hstr[y]\}]_m
$$
and thus \eqref{e-adj} follows. 
\end{pf}  
In other words $\hstr$ acts as a conjugation between multiplication and
contraction maps.

\subsection{Extension To Horizontal Forms}

The operator $\hstr$ will be extended to $\horPmc$ by left $\cal{B}$-linearity,
in other words we define
\begin{equation}
\hstrP(b\otimes\vartheta)=b\otimes\hstr(\vartheta) \qquad\quad\forall b\in\cal{B}\quad\forall
\vartheta\in\Vmc^\wedge. 
\end{equation}
Of course, here it is necessary to deal with extended horizontal forms $\horPmc$.
As a consequence of \eqref{v-hstr} we have
\begin{equation}\label{Fh-hstrP}
\Fh\hstrP=(\hstrP\otimes\id)\Fh. 
\end{equation}
This intertwining property, together with \eqref{circ-hstr} and the definition of the product 
in $\horPmc$, implies that $\hstrP$ is $\lambda$-twisted right
$\cal{B}$-linear, $\hstrP(\psi b_\lambda)=\hstrP(\psi)b$.

\begin{defn} The map $\hstrP$ is called the {\it Hodge $*$-operator} for $P$. 
\end{defn}

Let us observe that $\hstrP(\Wmc)=\Wmc$, as directly follows from 
\eqref{Fh-hstrP}. We shall denote by $\hstrM\colon\Wmc\rightarrow\Wmc$
the corresponding restricted map. 

The introduced integration map $\intP\colon\hor_P\rightarrow \Bbb{C}$
naturally extends, by left $\mc$-linearity, to $\intP\colon\horPmc\rightarrow \mc$. 
Such an extended map intertwines the actions of $G$ and satisfies 
$$
\intP[\varphi^*]=S\Bigl\{\intP[\varphi]^*\Bigr\},\qquad
\intP[\varphi q]=\intP[\varphi]S(q) \qquad\quad q\in \mc. 
$$

The Hodge $*$-operator, together with the extended integration map, enables us to introduce naturally
a scalar product in the algebra of horizontal forms. 

\begin{lem} The formula 
\begin{equation} 
<\varphi,\psi>=\intP \varphi^*\hstrP[\psi]
\end{equation}
defines a $\mc$-valued scalar product in $\horPmc$. This scalar product is
$G$-covariant, and in terms of the natural left $\cal{B}$-module identification
$\horPmc\leftrightarrow\cal{B}\otimes\Vmc^\wedge$ it is given by a direct
product of natural scalar products in $\Vmc^\wedge$ and $\cal{B}$.
\end{lem} 

\begin{pf} The $\Fh$-covariance of the introduced scalar product follows from 
\eqref{Fh-hstrP} and the $\Fh$-covariance of $\intP$. Let is further observe
that the defined scalar product possesses all appropriate $\mc$-anti/linearity
properties. This follows from left $\mc$-linearity of $\intP$ and $\hstr$ and
identities describing $S$-twisted right $\mc$-linearity of these maps.  

It is sufficient to verify the lemma on the non-correlated elements of the
extended horizontal forms algebra. A direct computation gives

\begin{multline*}
<b\otimes\vartheta,q\otimes\eta>=\intP\Bigl\{
(\vartheta^*\otimes b^*)(q\otimes\hstr[\eta])\Bigr\}
=\sum_{kl}\intP\Bigl\{b_k^*q_l\otimes \vartheta^*{\circ}(c_k^*
d_l)\hstr[\eta]\Bigr\}=\\
=\sum_{kl}\omega_P[b_k^*q_l] g_\wedge\bigl(\vartheta^*{\circ}(c_k^* d_l),
\eta\bigr)=\omega_P[b^*q]\langle\vartheta,\eta\rangle
\end{multline*} 
where $b,q\in\cal{B}$ and $\vartheta,\eta\in\Vmc^\wedge$. \end{pf}

In what follows we shall assume that $\horPmc$ is equipped with the constructed 
scalar product. The next lemma gives an explicit description of the (formal) adjoint 
covariant derivative map. 

\begin{lem} 
The map $\nabla\colon\horPmc\rightarrow\horPmc$ is adjointable, in other words
there exists a (necessarily unique) linear map $\nabla^\dagger\colon\horPmc
\rightarrow\horPmc$
such that 
\begin{equation}\label{adjnabla-def}
<\varphi,\nabla(\psi)>=<\nabla^\dagger(\varphi),\psi> \qquad\quad\forall\varphi,\psi\in\hor_P.
\end{equation}
Explicitly, 
\begin{equation}\label{nabla-adj}
\nabla^\dagger(\psi)=-\hstrP^{-1}\nabla\hstrP(\psi). 
\end{equation}
The map $\nabla^\dagger$ is $\mc$-bilinear. It intertwines the right action $\Fh$---in other words
\begin{equation}
\Fh\nabla^\dagger=(\nabla^\dagger\otimes\id)\Fh. 
\end{equation}
\end{lem}

\begin{pf} Let us start from the identity $\intP\nabla=0$. We compute 
\begin{multline*}
0=\intP\nabla\bigl(\varphi^*\hstrP(\psi)\bigr)\sim\intP\Bigl\{\nabla(\varphi)^*\hstrP(\psi)
\Bigr\}+\intP\Bigl\{\varphi^* \nabla\hstrP(\psi)\Bigr\}=\\=
<\nabla(\varphi),\psi>+<\varphi, \hstrP^{-1}\nabla\hstrP(\psi)>
\end{multline*}
and hence $\nabla$ is adjointable in $\hor_P$ and \eqref{adjnabla-def} holds. The 
$\Fh$-covariance property follows from the $\Fh$-covariance of $\nabla$ and the $\Fh$-invariance
of the scalar product. The $\mc$-linearity is a direct consequence of \eqref{nabla-adj}
and the $\mc$-linearity of $\nabla$. It also follows from \eqref{nabla-adj} and
the twisted $\mc$-linearity of $\hstr$. 
\end{pf}

There is a simple coordinate expression for the adjoint derivative. It is given by
\begin{equation}\label{nabla-adj-coor}
\nabla^\dagger=-\sum_{i=1}^d\partial_i \otimes \iota[\theta_i].
\end{equation}
This can be proved in several ways, for example it follows by explicitly
taking the adjoints of $\partial_i$ and $\theta_i$ in the coordinate
expression for $\nabla$. Let us observe that $\nabla^\dagger$ generically
takes the values from the extended horizontal algebra $\horPmc$. 

\subsection{Quantum Laplacian}
We are now ready to introduce a quantum Laplace operator, in a similar manner as in classical 
geometry. 

\begin{defn} A linear operator $\LapP\colon\horPmc\rightarrow\horPmc$ defined by 
\begin{equation}\label{LapP-def}
\LapP=\nabla\nabla^\dagger+\nabla^\dagger\nabla
\end{equation}
is called the {\it quantum Laplacian}. 
\end{defn}
By construction, $\LapP$ is a symmetric positive operator. Note that
it operates within the extended horizontal forms algebra. 
It is interesting to write down an explicit coordinate formula for
the quantum Laplacian.

\begin{pro} We have 
\begin{equation}\label{LapP-coor}
\LapP(b\otimes\vartheta)=-\sum_{ij}\partial_i \partial_j(b)\otimes g_{ij}\vartheta+
\frac{1}{2}\sum_{ij\alpha}b_\alpha \varrho_\nabla^{ij}(c_\alpha)\otimes
\theta_i\iota[\theta_j](\vartheta),
\end{equation}
where $F(b)=\Sum_\alpha b_\alpha\otimes c_\alpha$ and $g_{ij}=g(\theta_i,\theta_j)$.   
\end{pro} 

\begin{pf} Playing with coordinate formulas \eqref{nabla-coor} and 
\eqref{nabla-adj-coor}, and the commutation relations \eqref{com-rel-X}
and \eqref{contr-br} we find
\begin{equation*}
\begin{split}
\LapP(b\otimes\vartheta)&=-\sum_{ij}\Bigl\{\partial_i\partial_j(b)\otimes \iota[\theta_i] 
(\theta_j\vartheta)\Bigr\}-\sum_{ij}\Bigl\{\partial_i\partial_j(b)\otimes \theta_i\iota[\theta_j](\vartheta)\Bigr\}\\
&=-\sum_{ij}\partial_i\partial_j(b)\otimes g_{ij}\vartheta-\sum_{ij}\Bigl\{\partial_i\partial_j(b)\otimes 
\theta_i\iota[\theta_j](\vartheta)\Bigr\}\\{}&\phantom{=}{}+\sum_{ijkl}\partial_k\partial_l(b)\otimes\sigma_{kl}^{ij}
\theta_i\iota[\theta_j]
(\vartheta)\\
&=-\sum_{ij}\partial_i \partial_j(b)\otimes g_{ij}\vartheta+\frac{1}{2}\sum_{ij\alpha}
b_\alpha \varrho_\nabla^{ij}(c_\alpha)\otimes \theta_i\iota[\theta_j](\vartheta). \qed
\end{split}
\end{equation*}
\renewcommand{\qed}{}
\end{pf} 

In general, maps $g_{ij}\colon\Vmc^\wedge\rightarrow\Vmc^\wedge$ will be 
non-scalar operators. Here are some further elementary algebraic properties of 
$\LapP$---the covariance property 
\begin{equation}
\Fh\LapP=(\LapP\otimes\id)\Fh\label{Fh-LapP}, 
\end{equation}
and it is worth mentioning the following expressions
\begin{gather}
-\LapP\hstrP^{-1}=\nabla\hstrP\nabla+\hstrP^{-1}\nabla\hstrP\nabla\hstrP^{-1}\\
-\hstrP\LapP=\nabla\hstrP^{-1}\nabla+\hstrP\nabla\hstrP^{-1}\nabla\hstrP.
\end{gather}

According to \eqref{Fh-LapP} the map $\LapP$ is reduced in $\Omega_M$. We shall denote 
the corresponding restriction map by $\LapM\colon\Omega_M\rightarrow\Omega_M$. Obviously, 
\begin{gather*}
\LapM\dM=\dM\LapM\qquad \LapM\dM^\dagger=\dM^\dagger\LapM\\
\LapM=(\dM+\dM^\dagger)^2. 
\end{gather*}

Let us consider a natural decomposition
$$ \hor_P=\sideset{}{^\oplus}\sum_{\alpha\in\cal{T}}\cal{H}^\alpha 
\qquad\cal{H}^\alpha\leftrightarrow \hbim{\alpha}\otimes H_\alpha$$
into the multiple irreducible subspaces. It is worth mentioning that these subspaces 
are mutually orthogonal with respect to the scalar product $<,>$ in $\hor_P$. 
All the maps $\nabla,\nabla^\dagger$ and $\LapP$ are reduced in the spaces $\cal{H}^\alpha$. 
 
\section{Quantum Spin Bundles}

Throughout this section we shall assume that the structure group $G$ possesses a very special 
`spinorial' representation (and consequently we shall relax from the faithfulness assumption 
for $\v$). The corresponding frame structures on quantum spaces/bundles are then 
interpretable as `covering bundles' of the `real' orthonormal frame bundles. The orthogonal quantum group $G_0$ corresponds to the Hopf *-subalgebra $\cal{A}_0$ of $\cal{A}$ generated by the matrix elements $\v_{ij}$. 

The original orthonormal frame bundle $P_0$ is given by the *-subalgebra $\cal{B}_0$ of 
$\cal{B}$ generated by multiple irreducible submodules of $\cal{B}$ corresponding to 
the representations of $G_0$. Obviously $F(\cal{B}_0)\subseteq\cal{B}_0\otimes\cal{A}_0$ and 
$i(\cal{V})\subseteq\cal{B}_0$. Taking the corresponding restriction maps, 
we obtain a quantum principal $G_0$-bundle $P_0=(\cal{B}_0,i,F)$ over $M$ with the faithful 
action $\v$. Geometrically, $P$ is a kind of a covering space for $P_0$ and $P_0$ corresponds 
to the `vanilla' orthonormal frame bundle in classical geometry. 

We shall first formalize the idea of a `quantum spinor space'. We shall use the 
quantum Clifford algebra $\ClV$ associated to $\V$, metric coefficients algebra $\mc$, 
the braid operator $\sigma$ and quantum metric $g\colon\V\otimes\V\rightarrow \mc$. 
This algebra is constructed by $g$-deforming the product in $\Vmc^\wedge$,
while preserving the *-structure and the $\circ$-structure.

Let us assume that a finite-dimensional Hilbert space $\Salg$ is given, together with 
a unitary representation $\vspi\colon\Salg\rightarrow\Salg\otimes\cal{A}$. 
Let us also assume that $\Salg$ is an irreducible left *-module over $\ClV$. 
Finally, let us assume that the following compatibility condition holds:
\begin{equation}
\vspi(Z\xi)=\vmc(Z)\vspi(\xi)\qquad\qquad Z\in \ClV \quad\xi\in \Salg. 
\end{equation}
The meaning of this condition is that the action of $G$ on $\ClV$ can be viewed
as the adjoint action of $\vspi$, in terms of operators acting in $\Salg$. 
\begin{defn}
If the above conditions are fulfilled, we shall say that
$\Salg$ is a {\it quantum spinor space} associated to $G$ and $\ClV$.
\end{defn}

The $\ClV$-module structure $\gamma\colon \ClV\rightarrow B(\Salg)$ is
generally not faithful, and $\mc$ may be infinite-dimensional. 
The map $\gamma$ (including its values on $\mc$) is completely determined by 
the assignment $\gamma\colon\V\rightarrow B(\Salg)$. This simple observation
can be used as a starting point in constructing $\mc$ and $\ClV$. 

By combining the *-algebra structures on $\cal{B}$ and $\Vmc^\wedge\leftrightarrow\ClV$, we obtain a *-algebra $\ClP$. By  construction, 
we have a natural identification $\ClP\leftrightarrow \horPmc$ of $\cal{B}$-bimodules. 
The *-algebra structure on $\ClP$ is given by the standard cross-product type formulae
\begin{align*}
(q\otimes\vartheta)(b\otimes\eta)&=\sum_kqb_k\otimes(\vartheta\circ c_k)
\eta\\
(b\otimes\vartheta)^*&=\sum_k b_k^*\otimes(\vartheta^*\circ c_k^*). 
\end{align*}

By taking the product of the actions $\vcl$ and $F$ we obtain a *-homomorphism 
$\Fcl\colon\ClP\rightarrow\ClP\otimes\cal{A}$. Obviously $\Fcl\leftrightarrow\Fh$, 
in terms of the identification $$\ClP\leftrightarrow\horPmc.$$ We shall denote by 
$\ClM\subseteq\ClP$ the $\Fcl$-invariant *-subalgebra of $\ClP$. Obviously, we have 
a natural identification $$\ClM\leftrightarrow \Wmc.$$

Starting from a quantum spinor space, we can define the associated spinor bundle. 
Let us consider a free left $\cal{B}$-bimodule $\cal{S}$, given by
\begin{equation}
\cal{S}=\cal{B}\otimes \Salg. 
\end{equation}

By taking the product of actions $F$ and $\vspi$ we obtain the map 
$\Fspi\colon\cal{S}\rightarrow\cal{S}\otimes\cal{A}$. 
There exists a natural $\Fspi$-invariant scalar product on $\cal{S}$, defined by taking the 
direct product of the scalar products in $\Salg$ and $\cal{B}$. 
In what follows we shall assume that $\cal{S}$ is equipped with this scalar product. 
Let $\SpiM\subseteq\cal{S}$ be the subspace of $\Fspi$-invariant elements. This space is a 
$\cal{V}$-bimodule, in a natural way. In accordance with our general discussion,  it is interpretable as the appropriate associated spinor bundle.  

\begin{defn} The *-algebra $\ClM$ is called {\it quantum Clifford bundle algebra} over the space $M$. The $\cal{V}$-bimodule $\SpiM$ is called {\it quantum spinor bundle}, and its elements are called {\it quantum spinor fields} over $M$. 
\end{defn}

It is possible to introduce a natural action map $\beta\colon\ClP\otimes\cal{S}\rightarrow\cal{S}$ of 
$\ClP$ on $\cal{S}$. This map is defined by 
\begin{equation}
(q\otimes x)(b\otimes \zeta)=\sum_k qb_k \otimes (x{\circ} c_k)[\zeta]
\end{equation}
It is easy to see that this indeed defines a faithful unital action of $\ClP$ on $\cal{S}$, 
intertwining the corresponding natural coactions $\Fcl\times \Fspi$ and $\Fspi$. In particular, 
it follows that 
\begin{equation}
\beta(\ClM\otimes\SpiM)=\SpiM. 
\end{equation}

\begin{pro} \bla{i} The action $\beta\colon\ClP\otimes\cal{S}\rightarrow\cal{S}$ is hermitian, 
in other words 
\begin{equation}\label{ClP-adj}
<\psi,T\varphi>=<T^*\psi,\varphi>
\end{equation}
for each $\psi,\varphi\in\cal{S}$ and $T\in\ClP$. 

\smallskip
\bla{ii} The operators $T\colon\cal{S}\rightarrow\cal{S}$ coming from $\ClP$ are generally unbounded. However the restrictions on $\SpiM$ are bounded. In particular $\ClM$ acts on 
$\SpiM$ by bounded operators.  
\end{pro}

\begin{pf}
We compute
\begin{multline*}
<b\otimes \zeta, T(q\otimes \xi)>=\sum_\alpha 
<b\otimes\zeta, uq_\alpha \otimes [\vartheta{\circ}d_\alpha](\xi)>\\
=\sum_\alpha\omega_P(b^*uq_\alpha)(\zeta, [\vartheta{\circ}d_\alpha](\xi))
=\sum_\alpha\omega_P(b^*uq_\alpha)<[\vartheta{\circ}d_\alpha]^*(\zeta), \xi>\\
=\sum_\alpha\omega_P\bigl((u^*b)^*q_\alpha\bigr)<[\vartheta^*{\circ}\k(d_\alpha)^*](\zeta), \xi>\\=
\sum_{ik}\omega_P((u_i^*b_k)^*q)<[\vartheta^*{\circ}(c_k^*a_i)^*](\zeta), \xi>
=<T^*(b\otimes \zeta),q\otimes \xi>,
\end{multline*}
where $T=u\otimes\vartheta$ and the corresponding sums indicate the action of $F$ on 
the elements $q, b$ and $u$. Consequently \eqref{ClP-adj} holds. The fact that
the operators $T$ are unbounded in general, comes from the $\circ$-structure used in 
the definition of the action of $\ClP$. Indeed, the $\circ$ is generally {\it not continuous}
when viewed as a homomorphism ${\circ}\colon\cal{A}\rightarrow \mathrm{M}_d(\Bbb{C})$.  

On the other hand, multiplication operators by $q\in\cal{B}$ are continuous. If $\psi=\Sum_{b\zeta}b\otimes \zeta\in\SpiM$ then it is easy to see that 
$(q\otimes x)[\psi]=\Sum_{b\zeta}qb\otimes T_x[\zeta]$ which implies that 
the restricted action on $\SpiM$ is continuous. 
\end{pf}
 
\section{Quantum Dirac Operator}

Let us consider a linear operator $\D\colon\cal{S}\rightarrow\cal{S}$ given by 
\begin{equation}
\D(b\otimes x)=-i\sum_j\partial_j(b)\otimes \theta_j[x],
\end{equation}
where we have interpreted the elements of $\V$ as linear operators in $\Salg$, in 
accordance with our definition of quantum Clifford algebras. We begin by demonstrating a 
couple of elementary properties of the introduced map. 

\begin{pro} \bla{i} The map $\D$ is $\Fspi$-covariant. In other words,
\begin{equation}\label{D-Fspi}
\Fspi\D=(\D\otimes\id)\Fspi. 
\end{equation}
In particular $\D(\SpiM)\subseteq\SpiM$. 

\bla{ii} We have 
\begin{equation}\label{D-adj}
<\psi,\D(\varphi)>=<\D(\psi),\varphi>
\end{equation}
for each $\psi,\varphi\in\cal{S}$. 
\end{pro}

\begin{pf} The intertwining property of $\D$ is an immediate consequence of its
definition and the transformation properties of $\theta_i$ and $\partial_i$.  To check
\eqref{D-adj} let us observe that the scalar product in $\cal{S}$ is obtained by 
tensoring the scalar products in $\cal{B}$ and $\Salg$. Consequently, taking the formal adjoints, applying \eqref{X-adj} and playing with the basis $\theta_i$ we find
\begin{multline*}
\D^\dagger=i\sum_j(\partial_j\otimes\theta_j)^\dagger=
i\sum_j\partial_j^\dagger\otimes\theta_j^\dagger=\\=-i\sum_{jk}[C^{-1/2}_\v]_{kj}\partial_k
\otimes C^{-1/2}_\v(\theta_j)=-i\sum_k\partial_k\otimes\theta_k=\D. \qed
\end{multline*}
\renewcommand{\qed}{}
\end{pf}

\begin{defn} The map $\D\colon\SpiM\rightarrow\SpiM$ is called a {\it quantum Dirac operator}
for $M$. 
\end{defn}

The following simple proposition shows us that the Dirac operator contains the whole 
information about the differential $\dM\colon\Omega_M\rightarrow\Omega_M$, as in 
the classical theory. This fits into the axiomatic framework 
of \cite{c}. However, in contrast to \cite{c} the eigenvalues of our Dirac operator do not
obey the classical-type asymptotics in general. 

\begin{pro} We have 
\begin{equation}
[\D,f]=\dM(f)\qquad\quad\forall f\in\cal{V}. 
\end{equation}
Here we have interpreted $f$ and $\dM(f)$ as sections of the Clifford algebra bundle. 
\end{pro}

\begin{pf}
A direct computation gives 
\begin{multline*}
\D[f(b\otimes \vartheta)]=\sum_i\partial_i(fb)\otimes\theta_i\vartheta
=\Bigl\{\sum_i\partial_i(f)\otimes\theta_i\Bigr\}(b\otimes\vartheta)\\{}+
f\sum_i\partial_i(b)\otimes\theta_i\vartheta
=\dM(f)(b\otimes\vartheta)+f\D(b\otimes\vartheta). \qed 
\end{multline*}
\renewcommand{\qed}{}
\end{pf}

Let us now investigate the relation between $\D$ and the spinorial Laplacian. 

\begin{defn} In accordance with classical theory, the operator 
$\LapS\colon\SpiM\rightarrow\SpiM$ given by
\begin{equation}
\LapS=-\sum_{ij}\partial_i \partial_j\otimes g_{ij}
\end{equation}
is called the {\it spinorial Laplacian}. Here 
$g_{ij}=g(\theta_i,\theta_j)\colon\Salg\rightarrow\Salg$. 
\end{defn}
Here is a quantum version of classical Lichnerowicz formula. 
\begin{pro} We have 
\begin{equation}
\D^2=\LapS+ \widehat{R}
\end{equation}
where $\widehat{R}\colon\SpiM\rightarrow\SpiM$ is the `cliffordization' of the curvature. 
Explicitly $\widehat{R}\leftrightarrow\nabla^2$ in terms of the identification 
$\ClV\leftrightarrow\Vmc^\wedge$. 
\end{pro}

\begin{pf} This is a direct consequence of the definition of the quantum Dirac operator 
and the definition of the product in braided Clifford algebras. 
\end{pf}

\section{Examples $\&$ Remarks}

\subsection{General Connections}

So far we have only considered the torsionless Levi-Civita connections. Our results are 
easily applicable to more general connections. By definition, a {\it bundle derivative}
on a framed quantum principal bundle $P$ is an arbitrary first-order hermitian antiderivation
$D\colon\hor_P\rightarrow\hor_P$ extending the differential $\dM\colon\Omega_M\rightarrow
\Omega_M$ and intertwining the action map $\Fh$. The torsion tensor 
$T_D\colon\V\rightarrow\zh_P^2$ of a bundle derivative $D$ is defined by
$$ T_D(\theta)=D(\theta), $$
as in classical geometry. The curvature $\varrho_D\colon\cal{A}\rightarrow\zh_P^2$ of $D$ is 
given by $$ D^2(\varphi)=-(-)^{\partial\varphi}\sum_k\varphi_k\varrho_D(c_k).$$
The bundle derivatives form a real affine space $\der(P)$. 

\begin{lem} We have 
\begin{equation}
\intP D(\varphi)=0\qquad\qquad\forall D\in\der(P). 
\end{equation}
\end{lem}

\begin{pf} According to the results of \cite{d-frm2} we can write
$$ D(\varphi)=\nabla(\varphi)+(-)^{\partial\varphi}\sum_k\varphi_k \lambda(c_k), $$
where $\lambda\colon\cal{A}\rightarrow\zh_P$ is an appropriate linear map 
vanishing on scalars. Now
applying \eqref{intP-nabla} and performing elementary transformations we find
\begin{multline*}
\intP D(\varphi)=\intP \nabla(\varphi) +(-)^{\partial\varphi}\sum_k \intP [\varphi_k\lambda(c_k)]
=(-)^{\partial\varphi}\sum_k \intP [\varphi_k\lambda(c_k^{(1)})h(c_k^{(2)})]\\
=(-)^{\partial\varphi}\sum_k \intP [\varphi_kh(c_k)]\lambda(1)=0. \qed 
\end{multline*}
\renewcommand{\qed}{}
\end{pf}
 
In particular, this property implies that every bundle derivative will be formally adjointable, 
with the formal adjoint given by 
\begin{equation}
D^\dagger(\varphi)=-\hstrP^{-1} D\hstrP(\varphi). 
\end{equation}

Further proceeding with this, it is possible to write the explicit formulas for the
Dirac operator and the Laplacian associated to an arbitrary bundle derivative $D$. 
As explained in \cite{d-frm2} the system of bundle derivatives allows us to construct a
natural bicovariant *-calculus $\Gamma$ on the structure group $G$, and a natural calculus on 
the bundle (extending the algebra $\hor_P$ of horizontal forms by taking a cross product with 
the appropriate braided exterior algebra $\Gamma_{\inv}^\wedge$). We can also extend
the whole picture to the level of $\horPmc$ and $\Wmc$. 

\subsection{Quantum Hopf Fibration}

This highly instructive example shows us that the asymptotics of the eigenvalues of 
the Dirac operator could be quite surprising in the non-commutative context. We are going
to deal with the Dirac operator over a quantum 2-sphere. We refer to \cite{r-sfera} for detailed calculations. 

By definition, the quantum Hopf fibration is a quantum $\U(1)$-bundle over a quantum 2-sphere 
\cite{p}. The total space of the bundle is given by the quantum $\SU(2)$-group. 
In other words \cite{w-su2}, the bundle *-algebra $\cal{B}$  is
generated by two elements $\{\alpha,\gamma\}$ and the following relations:
\begin{gather*}
\alpha\alpha^*+\mu^2\gamma\gamma^*=1\qquad\alpha^*\alpha+\gamma^*\gamma=1\\
\alpha\gamma=\mu\gamma\alpha\qquad
\alpha\gamma^*=\mu\gamma^*\alpha\qquad
\gamma\gamma^*=\gamma^*\gamma,
\end{gather*}
where $\mu\in[-1,1]\setminus\{0\}$. On the other hand, the Hopf *-algebra $\cal{A}$ of the
structure group $G=\U(1)$ is generated by a single unitary element $U$. The coproduct is
specified by $\phi(U)=U\otimes U$.  

It is also worth mentioning that the matrix
\begin{equation*}
u=\begin{pmatrix}
\alpha&-\mu\gamma^*\\
\gamma&\alpha^*
\end{pmatrix}
\end{equation*} 
defines the fundamental representation of the quantum $\SU(2)$ group. 

The above relations defining $\cal{B}$ are equivalent to the unitarity property 
$u^{-1}=u^*$. The coproduct map $\phi\colon\cal{B}\rightarrow\cal{B}\otimes\cal{B}$
is uniquely determined by $\phi(u_{ij})=\Sum_k u_{ik}\otimes u_{kj}$. 
We have denoted by the same symbol for the coproducts on $G$ and $P$. 

Our structure group $G$ is understandable as a {\it subgroup} of $P$, 
in accordance with the identification 
$$\cal{A}\leftrightarrow\cal{B}/\mathrm{gen}(\gamma, \gamma^*).$$ If $\mu\neq 1,-1$ then
$G$ is exactly the classical part of $P$. The action $F\colon\cal{B}\rightarrow\cal{B}\otimes\cal{A}$ 
is constructed from the coproduct, by taking the factor-projection on the second tensoriand, in
other words $F=(\id\otimes [\,]_{\cal{A}})\phi$. The algebra $\cal{V}$ describing the quantum
2-sphere $M$ is defined as the $F$-fixed point subalgebra of $\cal{B}$, so that the map $i$ is just 
the inclusion. 

Let us now sketch the construction of a canonical frame structure on the quantum 2-sphere 
\cite{d-affi}. We shall start from the canonical $3$-dimensional left-covariant and *-covariant 
calculus $\Phi$ over $P$. This calculus is constructed in \cite{w-su2}.
Explicitly the space $\Phi_{\inv}$ is spanned by the elements 
$$\eta_3=\pi(\alpha-\alpha^*)\qquad
\eta_+=\pi(\gamma)\qquad
\eta_-=\pi(\gamma^*)$$
and the canonical right ${\cal B}$-module structure $\circ$ on $\Phi_{\inv}$ is given by 
\begin{equation*} 
\begin{aligned}
\mu^2\eta_3{\circ}\alpha&=\eta_3\\
\mu\eta_\pm{\circ}\alpha&=\eta_\pm
\end{aligned}\qquad
\begin{aligned}
\eta_3{\circ}\alpha^*&=\mu^2\eta_3\\
\eta_\pm{\circ}\alpha^*&=\mu\eta_\pm
\end{aligned}\qquad\Phi_{\inv}{\circ}\gamma=\Phi_{\inv}{\circ}\gamma^*=\{0\}.
\end{equation*}

This $\cal{B}$-module structure on $\Phi_{\inv}$ factorizes through the ideal 
$\mathrm{gen}\{\gamma, \gamma^*\}$ and induces a right $\cal{A}$-module structure on 
the same space (and will be denoted by the same symbol). Let $\V$ be a vector space spanned 
by $\eta_\pm$. It will be equipped with the constructed $\circ$ and *-structures, and 
we shall assume that 
\begin{equation}
\chi(\eta_+)=\eta_+\otimes U^2\qquad\quad\chi(\eta_-)=\eta_-\otimes U^{-2}.
\end{equation}
Such a definition allows us to interpret $\chi\colon\V\rightarrow\V\otimes \cal{A}$ 
as {\it the adjoint} action of $G$, coming from the group structure in $P$.  
It follows that (in the basis $\eta_\pm$) the associated 
braid operator $\tau\colon\V^{\otimes 2}\rightarrow\V^{\otimes 2}$ looks like
\begin{equation}
\tau=\begin{pmatrix}1/\mu^2& 0 & 0 & 0\\
                         0 & 0 & 1/\mu^2 & 0\\
                         0 & \mu^2 & 0 & 0\\
                         0 & 0 & 0 & \mu^2 \end{pmatrix}
\end{equation}
The corresponding $\tau$-exterior algebra is given by the relations 
\begin{equation}
\eta_\pm^2=0\qquad\quad\eta_+\eta_-=-\mu^2\eta_-\eta_+. 
\end{equation}

It is worth noticing that these relations are a subset of the relations defining 
the canonical higher-order calculus over $P$, given by the universal differential
envelope $\Phi^\wedge$ of $\Phi$ (\cite{w2},\cite{d1}--Appendix B). For completeness, 
we shall list here the remaining relations (involving $\eta_3$). These are
\begin{equation}
\eta_3^2=0, \qquad\quad \eta_3\eta_\pm=\mu^{\mp 4}\eta_\pm\eta_3.
\end{equation} 
This means that $\hor_P$ is viewable as a subalgebra of $\Phi^\wedge$ generated by 
$\cal{B}=\hor_P^0$ and the elements $\eta_+,\eta_-$. Furthermore, we can 
introduce a natural projection homomorphism 
$p_{h\!or}\colon \Phi^\wedge\rightarrow \hor_P$, defined by 
$$ p_{h\!or}(\eta_3)=0\qquad p_{h\!or}\restr\hor_P=\id. $$
The canonical antiderivation $\nabla\colon\hor_P\rightarrow\hor_P$ is defined as 
the composition of this projection with the differential 
$d\colon\Phi^\wedge\rightarrow\Phi^\wedge$. 

The constructed map coincides with the covariant derivative of the canonical 
regular connection introduced in \cite{d2}. It corresponds to the standard 
Levi-Civita connection on the 2-sphere. 

Let us observe that we are already in the context of the spin bundles. The analog of
the orthonormal frame bundle $P_0$ over the 2-sphere $M$ is given by the *-subalgebra 
$\cal{B}_+\subseteq\cal{B}$ corresponding to the quantum $\SO(3)$ group (even 
combinations of generators $\alpha, \gamma,\gamma^*\alpha^*$). The structure group 
$G=\U(1)$ is here understood as a 2-fold covering of the structure group $G_0=\SO(2)$ of $P_0$. 

A braid operator $\sigma$ is given by the matrix
\begin{equation}
\sigma=\begin{pmatrix}1 & 0 & 0 & 0\\
                         0 & 0 & 1/\mu^2 & 0\\
                         0 & \mu^2 & 0 & 0\\
                         0 & 0 & 0 & 1 \end{pmatrix}
\end{equation}
and the *-structure on $\V$ is specified by $\eta_+^*=\mu\eta_-$, while
$\eta_3^*=-\eta_3$. The quantum metric is defined by 
\begin{equation}
g_{+-}^*=g_{+-} =\mu^2 g_{-+}\qquad\qquad g_{++}=g_{--}=0. 
\end{equation}
The algebra $\mc$ is infinite-dimensional. A realization of $\mc$ in the Hilbert
space $H=l^2(\Bbb{Z})$ is given by 
\begin{equation}
g_{+-}:e_k\mapsto\frac{1}{2}\mu^{2k}e_k
\end{equation}
where $\{e_k\bigm\vert k\in\Bbb{Z}\}$ are canonical basis vectors in $H$. A common domain 
$H_0$ for all the operators from $\mc$ consists of sequences with finite support, 
for example. 

In what follows we shall assume that the deformation parameter is positive (the 
case $\mu <0$ gives the same results). 
For the spinor space, we shall take $\Salg=\Bbb{C}^2$
with the canonical basis $|{+}\rangle,|{-}\rangle$ and the action
$$
\vspi|{+}\rangle=|{+}\rangle\otimes U\qquad\quad
\vspi|{-}\rangle=|{-}\rangle\otimes U^{-1}.
$$
The metric $g$ and
the algebra $\mc$ are completely determined by the assignment
\begin{equation}
\eta_+ \longmapsto \mu^{1/2}\begin{pmatrix}0 & 1\\0 & 0 \end{pmatrix}\qquad
\eta_- \longmapsto \mu^{-1/2}\begin{pmatrix}0 & 0\\1 & 0 \end{pmatrix}
\end{equation}
in particular it follows that 
\begin{equation}
\frac{1}{\mu}\gamma[g_{+-}]=\mu \gamma[g_{-+}]
=\frac{1}{2}\begin{pmatrix}\mu^{-1}& 0\\
0& \mu\end{pmatrix}
\qquad \gamma[\mc]=\Bbb{C}\oplus\Bbb{C}. 
\end{equation}
We see that the spinor representation $\gamma\colon\mc\rightarrow L(\Salg)$ is not faithful. 

By the way, it is easy to see that all possible irreducible representations of 
this Clifford algebra (with fixed spectral properties of $\mc$) are given by
$$
\gamma[\eta_+]= \mu^{1/2+k}\begin{pmatrix}0 & 1\\0 & 0 \end{pmatrix}\qquad
\gamma[\eta_-]= \mu^{-1/2+k}\begin{pmatrix}0 & 0\\1 & 0 \end{pmatrix}
$$
where $k\in\Bbb{Z}$. 

We are now going to study the associated Dirac operator.
For each spin level $s\in\Bbb{N}/2$ let $u^s$ be the matrix of the canonical spin-$s$ 
irreducible representation of $P$.  The matrix elements of all possible representations 
$u^s$ form a natural basis in $\cal{B}$. Let us denote by $\cal{B}_s$ the subspace of 
$\cal{B}$ spanned by the matrix elements of $u^s$. We have 
$$ \cal{B}=\Bbb{C}\oplus \sideset{}{^\oplus}\sum_{s\in \Bbb{N}/2}\cal{B}_s.$$ 

By construction, the coordinate vector fields $\partial_-,\partial_+$ of the frame structure coincide
with the spin creation and annihilation operators $i K_\mp$ for the right regular representation, 
given by the coproduct $\phi\colon\cal{B}\rightarrow\cal{B}\otimes\cal{B}$. Obviously 
$\phi(\cal{B}_s)\subseteq\cal{B}_s\otimes\cal{B}_s$, 
and moreover the space $\cal{B}_s$ is characterized as the multiple irreducible spin-$s$ 
subspace of $\cal{B}$. In terms of the matrix elements $u^s_{ij}$, 
the operators $\partial_\pm$ act nontrivially only on the second indexes, while the 
first indexes are `free'. 

Therefore, we can write 
$$\cal{B}_s\leftrightarrow \underbrace{H_s\oplus\dots\oplus H_s}_{2s+1}$$
and introduce a canonical basis in $\cal{B}_s$ of the form 
$\bigl\{\psi_{\alpha s}^m\bigm\vert m,\alpha=-s,\dots,s\bigr\}$. Here $\alpha$ is interpreted as a `degeneration index'. 

In summary, we have 
\begin{align*}
K_+(\psi_{\alpha s}^m)=-i\partial_+(\psi_{\alpha s}^m)&=\frac{1}{\mu^{s+m}}\Bigl\{(s-m)_\mu(s+m+1)_\mu\Bigr\}^{1/2}\psi_{\alpha s}^{m+1}\\
K_-(\psi_{\alpha s}^m)=-i\partial_-(\psi_{\alpha s}^m)&=\frac{\mu}{\mu^{s+m}}\Bigl\{(s-m+1)_\mu(s+m)_\mu\Bigr\}^{1/2}\psi_{\alpha s}^{m-1}\\
F(\psi_{\alpha s}^m)&=\psi_{\alpha s}^m\otimes U^m \qquad\quad
n_\mu=\frac{1-\mu^{2n}}{1-\mu^{2\phantom{n}}}. 
\end{align*}
Hence the spinor module $\cal{S}$ is decomposed as follows 
\begin{equation}
\cal{S}=\sideset{}{^\oplus}\sum_{s\in\Bbb{N}-1/2} \cal{S}_s, 
\end{equation}
where the spaces $\cal{S}_s$ are spanned by vectors
$$ 
\psi_{\alpha s}^{1/2}\otimes |{-}\rangle\qquad\quad \psi_{\alpha s}^{-1/2}\otimes |{+}\rangle
$$
with the degeneracy index $\alpha$ arbitrary. 

The Dirac operator is given by 
$$
i\D=\partial_+\otimes \begin{pmatrix} 0 & 0\\ 1 & 0\end{pmatrix}+\partial_-\otimes\begin{pmatrix}
0 & 1\\ 0 & 0\end{pmatrix}
$$
and it follows that $\D(\cal{S}_s)\subseteq\cal{S}_s$ for each $s\in\Bbb{N}-1/2$. It is now
very easy to diagonalize the reduced operators. We have two eigenvalues with 
the eigenvectors of the form
\begin{equation}
\D\varphi_{\alpha s}^+=\lambda_s\varphi_{\alpha s}^+\qquad\qquad
\D\varphi_{\alpha s}^-=-\lambda_s\varphi_{\alpha s}^-,
\end{equation}
where 
\begin{equation}
\lambda_s=\frac{\mu^{s+1/2}-\mu^{-s-1/2}}{\mu-\mu^{-1}}.
\end{equation}
It is worth noticing
that the eigenvalues satisfy the following recurrent formula
\begin{equation}
\lambda_s=\frac{\lambda_{s+1}+\lambda_{s-1}}{\mu+\mu^{-1}}. 
\end{equation}
In the above formulas we assumed that $\mu>0$. The spectrum is invariant if we 
replace $\mu\mapsto -\mu$. On the other hand, if $\mu=1,-1$ the spectrum will be 
given by $\lambda_s=2s$. The case $\mu=-1$ is very special, as it gives us a 
quantum spin bundle over the classical 2-sphere $M$. It illustrates a general phenomenon
that the classification problem of quantum principal bundles is qualitatively different from
its classical counterpart, even if both the base manifold and the structure group are
classical. 

\subsection{Discussion of Quantum Phenomena}

In classical geometry a very important geometrical information is contained in the 
eigenvalues of standard elliptic operators. In particular, the distribution of eigenvalues of 
the Dirac operator reflects the deepest structure of compact Riemannian spin manifolds. 

If we look at the asymptotics of the spectrum of (the modulus of) the Dirac operator over 
the quantum 2-sphere (with their degeneracies taken into account) we arrive 
to the following expression
\begin{equation}
a_N\sim 
\begin{cases} \sqrt{N} & \text{for $\mu=\pm 1$,}\\
\mu^{-\sqrt{N/2}} & \text{if $\mu\in(0,1)$}
\end{cases}
\end{equation}
with the symmetry $\mu\mapsto -\mu$. In particular, we see that the inverse of the 
Dirac operator is {\it trace class} in the fully quantum 
case $\mu\neq 1,-1$. This is not compatible with the formulation proposed in \cite{c}, 
where it was assumed that the quantum Dirac operator will always have a similar 
asymptotics as in the classical geometry 
$$a_N\sim N^{1/d}, \qquad\quad d=\dim(M).$$ 

Quantum geometry gives us much more freedom, and it is not possible to cover the diversity 
of all possible quantum spaces by a single asymptotic expression. 

In our theory (as far as we consider the pure Riemannian geometry) the group $G$ plays the 
role of special orthogonal structure group
$\SO(d)$, and $\v$ plays the role of its fundamental representation. 
In various interesting examples the representation $\v\colon\V\rightarrow\V\otimes\cal{A}$ will be 
irreducible. However, in general this representation will be reducible.
The fact that $\mc\neq \Bbb{C}$ allows us to overcome inherent obstacles that
would appear in the formalism, in the case of irreducible $\v\colon\V\rightarrow\V\otimes\cal{A}$. 

This is due to the fact that the braid operators $\tau,\sigma\colon\V\otimes\V\rightarrow\V\otimes\V$, the
*-structure and the quadratic form on $\V$ are all $\v$-covariant. If the metric
$g$ would take values from $\Bbb{C}$ then using elementary algebraic
operations with $g$, $*$ and $\sigma,\tau$ we would be able to built 
various intertertwiners of $\v$ out of these objects.
In general, these intertwiners would be non-scalar operators,
which implies that $\v$ will be reducible. The non-triviality of $\mc$
overcomes this obstacle. 

Another interesting quantum phenomenon is that the grade $m$ of the volume form $w$ 
is not necessarily the same as the number $d$ of coordinate one-forms. 

Our definition of a quantum Clifford algebra is motivated by considerations presented in
\cite{ozi} and \cite{d-ozi}, based on deformations of braided exterior algebras. Our main 
condition is similar---the vector space $\V^\wedge$ is equipped with a new product 
however this new product is a quantization of $\wedge$ with a non-commutative deformation
parameter. Of course, this is connected with a non-scalar nature of metric
components $g_{ij}\colon\Salg\rightarrow\Salg$. 

Finally, let us mention that among possible algebraic conditions involving braid 
operators $\tau$ and $\sigma$, there exist natural equations \cite{d-elliptic}
closely related to multiplicative unitaries \cite{w-multu} associated to 
compact quantum groups. 

\appendix
 
\section{Simple C*-algebraic Considerations}

In this appendix we shall consider the vertical integration map, and the associated GNS-type 
construction. Let $P=(\cal{B},i,F)$ be an arbitrary quantum principal $G$-bundle over $M$. 

Let us consider the map $\vint\colon\cal{B}\rightarrow\cal{V}$ defined by 
$i[\vint(b)]=(\id\otimes h)F(b)$. We shall assume that $\cal{V}$ is realized as an everywhere
dense *-subalgebra of a unital C*-algebra $\cstV$. 

One of the principal properties we have considered in the main text was the strict positivity 
condition. We shall first prove that this condition will be automatically satisfied in certain 
important special cases. 

\begin{pro} Let us assume that $\cal{B}$ admits a C*-algebraic closure $\cstB$ such that 
$i\colon\cal{V}\rightarrow\cal{B}$ extends to a *-monomorphism $i\colon\cstV\rightarrow\cstB$. 
Let us assume that $\cal{V}$ is stable in $\cstV$, under holomorphic functional
calculus. Then the vertical integration map $\vint$ is strictly positive. 
\end{pro}

\begin{pf} Let us start from the canonical \cite{w2} orthogonality relations
\begin{align}
h[u_{ki}^*u_{lj}]&=\delta_{ij}[C^{-1}_u]_{lk}/\tr(C_u)\\
h[u_{ki}u_{lj}^*]&=\delta_{kl}[C_u]_{ji}/\tr(C^{-1}_u)
\end{align}
where $u$ is an arbitrary unitary irreducible matrix representation of $G$, and $C_u$ is the
canonical intertwiner between $u$ and its second contragradient $u^{c c}$. 

In order to verify the positivity property, it is sufficient
to consider only the elements from a fixed multiple irreducible submodule $\cal{B}^u$.

A general element of such a form is given by
$$ 
b=\sum_{i\alpha} f_{i \alpha}b_{\alpha i}=\tr[\Phi B]
$$ 
where $f_{i\alpha}\in\cal{V}$ form a matrix $\Phi$ while $b_{\alpha i}$ form a matrix $B$
and satisfy
\begin{equation}\label{B-mult}
\sum_\alpha b_{\alpha i}^* b_{\alpha j}=\delta_{ij}1 
\qquad F(b_{\alpha i})=\sum_j b_{\alpha j}\otimes u_{ji}.
\end{equation}
This is a consequence of the holomorphic stability of $\cal{V}$ in $\cstV$. 
Now a direct computation gives
\begin{multline*}
\vint (b b^*)=\sum_{ij\alpha\beta} f_{i\alpha}\vint(b_{\alpha i}b_{\beta j}^*) f_{j\beta}^*
=\sum_{ijpq\alpha\beta} f_{i\alpha}b_{\alpha p} b_{\beta q}^*h(u_{pi}u_{qj}^*)f_{j\beta}^* \\
=\sum_{ijpq\alpha\beta} f_{i\alpha}b_{\alpha p} b_{\beta q}^*f_{j\beta}^*\delta_{pq}[C_u]_{ji}/ \tr[C_u^{-1}]\\=\sum_{ij\alpha\beta} f_{i\alpha}(B B^\dagger)_{\alpha\beta}f_{j\beta}^*
[C_u]_{ji}/\tr[C_u^{-1}]=\tr\Bigl\{C_u\Phi B B^\dagger\Phi^\dagger\Bigr\}/\tr[C_u^{-1}]
\end{multline*}
and from this expression it follows that $\vint$ is really strictly positive. 
\end{pf}

Let us now assume that the bundle $P$ is such that $\vint$ is strictly positive. 
\begin{pro} This property of $\vint$ implies the existence of 
a canonical faithful GNS-type representation of $\cal{B}$, which enables us to introduce a 
C*-norm in $\cal{B}$. 
\end{pro}

\begin{pf}
The strict positivity property implies that we can introduce a pre-Hilbert C*-module 
$\cal{B}$, with the $\cal{V}$-valued scalar product 
$$ <b,q>_{\! M}=\vint(b^* q). $$
The multiple irreducible subspaces $\cal{B}^u$ will be orthogonal 
relative to $<,>_{\! M}$. Let us denote by $\cal{H}$ the completion of this module. The 
multiplication map in $\cal{B}$ naturally induces a *-homomorphism 
$D\colon\cal{B}\rightarrow\Bbb{B}(\cal{H})$. Here $\Bbb{B}(\cal{H})$ is the C*-algebra of bounded right $\cstV$-linear adjointable operators in $\cal{H}$. In order to 
prove the boundness of the operators $D[\,]$ it is sufficient to check it on canonical 
multiplets $b_{\alpha i}$ and on the elements $f\in \cal{V}$. We have 
\begin{multline*}
<D[b_{\alpha i}]\psi, D[b_{\alpha i}]\psi>_{\! M}\leq \sum_\beta <D[b_{\beta i}]\psi, 
D[b_{\beta i}]\psi>_{\! M}=\sum_\beta <\psi, D[b_{\beta i}^*b_{\beta i}]\psi>_{\! M}\\
=<\psi,\psi>_{\! M}
\end{multline*}
and it follows that $D[b_{\alpha i}]$ are bounded with $|D[b_{\alpha i}]|\leq 1$. It is 
easy to verify that $D\colon \cal{V}\rightarrow \Bbb{B}(\cal{H})$ is isometric. 

Obviously, the unit element $1\in\cal{H}$ is a cyclic and separating vector for $D$. In particular $D$ is faithful. 
\end{pf}

The constructed C*-norm has some further interesting properties. At first, let 
us observe that 
\begin{equation}
\bigm\vert \vint b\bigm\vert\leq \bigm\vert b\bigm\vert 
\end{equation}
for each $b\in\cal{B}$. In particular, the vertical integration map extends continuously 
to the C*-completion $\cstB$ and the above inequality holds on the whole $\cstB$. Secondly, 
the C*-structure on $\cal{B}$ is unique. Indeed, let $|\,|_{\mathrm{max}}$ be the 
universal C*-norm on $\cal{B}$ extending the C*-norm on $\cal{V}$, 
and let $B_{\mathrm{max}}$ be the corresponding C*-completion. The existence of such a maximal
C*-norm follows easily from equalities \eqref{B-mult} and the existence of $D$. By 
construction $F$ extends by continuity to a *-homomorphism $F\colon 
B_{\mathrm{max}}\rightarrow B_{\mathrm{max}}\otimes A$. It is easy to see that $\omega=(\id\otimes h)F\colon B_{\mathrm{max}}\rightarrow \cstV$ is strictly positive 
(assuming that the Haar measure on $G$ is faithful) and 
$|\omega(b)|\leq |b|_{\mathrm{max}}$ on $B_{\mathrm{max}}$. This implies that $B_{\mathrm{max}}=\cstB$. In particular $|\,|_{\mathrm{max}}=|\,|$ and $\omega=\vint$. 

The map $F\colon \cstB\rightarrow\cstB\otimes A$ satisfies 
\begin{equation*}
\begin{CD}
\cstB @>{\mbox{$F$}}>> \cstB\otimes A\\
@V{\mbox{$F$}}VV     @VV{\mbox{$\id\otimes \phi$}}V\\
\cstB\otimes A @>>{\mbox{$F\otimes\id$}}> \cstB\otimes A\otimes A
\end{CD}\qquad\qquad
\cstB\otimes A=\lin\overline{\Bigl\{F(b)a\Bigr\}} 
\end{equation*}
which are purely C*-algebraic counterparts of the quantum group action axioms. The 
second equality expresses the idea that the group $G$ acts by homeomorphisms on the bundle 
space $P$. We conclude this appendix by 
\begin{equation*}
\cstB\otimes A=\lin\overline{\Bigl\{qF(b)\Bigm\vert q,b\in \cstB\Bigr\}} 
\end{equation*}
which is a purely C*-algebraic version of the freeness axiom. 

\section{Braided Clifford Algebras}

In this Appendix we present an introduction to a general
theory of braided Clifford algebras and spinors. Much of the material
is logically independent of our main context of principal bundles,
quantum groups and bicovariant bimodules.

Let us consider a complex finite-dimensional vector space $\V$ equipped
with a regular braid operator $\sigma\colon\V\otimes\V\rightarrow\V\otimes\V$ and 
an antilinear involution $*\colon\V\rightarrow\V$. The involution map extends
naturally to a *-structure on the tensor algebra $\V^\otimes$. 
Explicitly,
$$ (\xi_1\otimes\dots\otimes\xi_n)^*=\xi_n^*\otimes\dots\otimes\xi_1^*. $$
We have denoted by the same symbol $*$ a unique
antimultiplicative (unital) antilinear extension on the tensor
algebra. We shall assume that
\begin{equation}\label{*-sigma} *\sigma=\sigma *. 
\end{equation}

First, we are going to formalize the idea of a {\it quantum metric}. As
we already mentioned it will be allowed that `metric coefficients'
do not commute.

Let $\mc$ be a *-algebra, and let $g\colon\V\otimes\V\rightarrow\mc$ be
a linear map. Below we have listed a number of interesting identities
involving $g$, $*$ and $\sigma$: 

\smallskip
\bla{i} {\bf Braided-symmetricity of the metric.} In other words
\begin{equation}
g\sigma=g.
\end{equation}
For this to make any sense, it is necessary that $1$ belongs to the
spectrum of $\sigma$. 

\bla{ii} {\bf Reality property.} We have
\begin{equation}
g(x,y)^*=g(y^*,x^*)\qquad \forall x,y\in\V
\end{equation}

\bla{iii} {\bf Funny $\sigma$-compatibility.} We have
\begin{equation}\label{g-g}
\begin{aligned}
g\otimes_\mc g&=(g\otimes_\mc g)(\id\otimes\sigma\otimes\id)(\sigma^{-1}\otimes\sigma)
(\id\otimes\sigma^{-1}\otimes\id)\\
g\otimes_\mc g&=(g\otimes_\mc g)(\id\otimes\sigma\otimes\id)(\sigma\otimes\sigma^{-1})
(\id\otimes\sigma^{-1}\otimes\id).
\end{aligned}
\end{equation}
Let us observe that the above two equations are {\it equivalent} if we assume 
that reality condition holds. As we shall see, this property ensures that $g$
is extendible, to the level of appropriate $\mc$-bimodules. 

\bla{iv} {\bf Weak positivity.} In order to formulate this property we have to assume
that $\mc$ is realized by operators in the Hilbert space $H=l^2(\Bbb{Z})$. Since in 
general these operators will be {\it unbounded} we have to take care about the domains. 
We shall assume that there is an everywhere dense linear subspace $H_0\subseteq H$ which is 
a common domain for all the operators from $\mc$. We shall also assume that the *-structure
on $\mc$ is represented as taking formal adjoints of linear operators in $H_0$, 
and that there exists cyclic and separating vectors $\Omega\in H_0$ for $\mc$. 

So one natural version of positivity would be
\begin{equation}
g(x^*,x)\geq 0\qquad x=0 \Leftrightarrow g(x^*,x)=0\qquad\forall x\in\V.   
\end{equation}
The reason why we call this condition `weak' positivity will become
clear after we construct a canonical $\mc$-bimodule structure over $\V$,
and introduce a stronger version of positivity. 

\smallskip
\bla{v} {\bf Minimality $\&$ Invertibility.} The matrix $g_{ij}=g(\theta_i,\theta_j)$ is 
invertible in $\mc$, where $\{\theta_1,\dots,\theta_d\}$ are basis vectors in $\V$. 
Moreover, the algebra $\mc$ is generated by the matrix elements of $g$ and $g^{-1}$.

\smallskip
\bla{vi} {\bf Twisting $\mc$ and $\V$.} Let $\mcVt\colon\mc\rightarrow\mathrm{M}_d(\mc)$ 
be a unital homomorphism. Obviously, this map is completely
determined by its values on the elements $g_{ij}$ and it gives us the structure of 
a right $\mc$-module, in the free left $\mc$-module $\Vmc\leftrightarrow\mc\otimes\V$,
so that we have a bimodule structure. The right $\mc$-multiplication is simply given by 
$$ \theta_i q=\sum_j \mcVt(q)_{ij}\theta_j. $$
Equivalently, the map $\mcVt$ can be viewed as a twisting operator 
$$ \mcVt\colon \V\otimes\mc\rightarrow\mc\otimes\V \qquad \mcVt(\theta_i\otimes q)=\Sum_j\mcVt(q)_{ij}\otimes \theta_j. $$
This twisting preserves the product and the unit in $\mc$, in a natural way. 

We shall assume that 
\begin{equation}
\mcVt(\id\otimes g)=(g\otimes \id)(\id\otimes\sigma)(\sigma^{-1}\otimes\id), 
\end{equation}
which completely fixes $\mcVt$. 

\begin{lem} The *-involutions on $\V$ and $\mc$ naturally combine to a 
*-structure on the bimodule $\Vmc$. In particular the map $\mcVt\colon \V\otimes\mc\rightarrow\mc\otimes\V$ is invertible and 
\begin{equation}\label{mcVt-*}
*\mcVt*=\mcVt^{-1}. 
\end{equation} 
In particular, it follows that $\Vmc$ free, as a right $\mc$-module. 
\end{lem}

\begin{pf} The map $*\colon\Vmc\rightarrow\Vmc$ is introduced by 
$(q \otimes \theta_j)^*=\theta_j^*q^*$. It is sufficient to prove that such a map is 
involutive. This follows easily from \eqref{mcVt-*} and the reality properties for 
$\sigma$ and $g$. 
\end{pf}

More generally, it is easy to see that 
$$
\Vmc^{\otimes n}\leftrightarrow \overbrace{\V\otimes\dots\otimes\V}^n\otimes\mc
\leftrightarrow \mc\otimes\underbrace{\V\otimes\dots\otimes\V}_n
$$
in a natural manner, as a right/left $\mc$-module. 

Let us observe that our definition of the $\mc$-bimodule structure implies
that the braiding $\sigma\colon\V\otimes\V\rightarrow\V\otimes\V$ is (necessarily uniquely) 
extendible to a bimodule homomorphism 
$\sigma\colon\Vmc^{\otimes 2}\rightarrow \Vmc^{\otimes 2}$. Similarly, property 
\eqref{g-g} ensures that the metric tensor $g\colon\V\otimes\V\rightarrow\mc$ is 
uniquely extendible to a $\mc$-bilinear map $g\colon\Vmc^{\otimes 2}\rightarrow\mc$. 
In playing with such extended maps, it is useful to recall the following simple 
lemma:

\begin{lem} Let $\cal{M}$ and $\cal{N}$ be left $\mc$-modules and let 
$\Phi\colon\cal{M}\rightarrow\cal{N}$ be a linear map satisfying 
$$
\Phi(g_{ij}\xi)=g_{ij}\Phi(\xi). 
$$
Then $\Phi$ is left $\mc$-linear. 
\end{lem}

\begin{pf} It is sufficient to check that $\Phi$ commutes with left multiplications by
the matrix elements $[g^{-1}]_{ij}$. However this is equivalent to the above formula. 
\end{pf}

Our extensions preserve all relevant algebraic relations between $g$, $\sigma$ and $*$. 
There is an interesting way to describe the relation between 
the left and the right $\mc$-module structures on $\Vmc$. They are related by
\begin{equation}\label{funny}
(\id\otimes g)(\sigma\otimes\id)=(g\otimes\id)(\id\otimes\sigma).
\end{equation}
In this formula the tensor product is taken over $\mc$. 

The space $\Vmc$, together with the extended $\sigma$, generates a braided monoidal 
category $\cal{C}$. We shall use the same symbol $\sigma$ to denote the generic 
braiding in this category. Moreover, we shall use the same symbol $g$ for extended
contraction maps $g\colon\Vmc^{\otimes n}\otimes_{\mc}\Vmc^{\otimes n}\rightarrow\mc$. 
The extended maps are defined inductively by
\begin{equation}
g\bigl\{(\psi\otimes x)\otimes (y\otimes\xi)\bigr\}=g(\psi,g(x,y)\xi), 
\end{equation}
where $x,y\in\V$ and we shall also assume that tensors with different grades are mutually 
`orthogonal'. 

Let us now consider a map $\langle\rangle:\Vmc\times\Vmc\rightarrow\mc$ defined by: 
\begin{equation}
\langle\psi,\xi\rangle=g(\psi^*,\xi). 
\end{equation}
It is easy to see that the following identities are fulfilled: 
\begin{gather}
\langle\psi,\xi a\rangle=\langle\psi,\xi\rangle a \qquad\langle\psi a,\xi\rangle=a^*\langle\psi,\xi\rangle\\
 \langle\psi,a\xi\rangle=\langle a^*\psi,\xi\rangle\\
\langle\psi,\xi\rangle^*=\langle\xi,\psi\rangle\\
\langle\psi,\xi+\varphi\rangle=\langle\psi,\xi\rangle+\langle\psi,\varphi\rangle. 
\end{gather} 
This map plays the role of a hermitian scalar product in $\Vmc$. To completely justify
this interpretation, it is necessary to formulate the appropriate positivity condition. 

\smallskip
\bla{vii} {\bf Strict positivity.} Assuming that $\mc$ is realized in $H=l^2(\Bbb{Z})$ 
we have 
\begin{equation}
\langle\xi,\xi\rangle=g(\xi^*,\xi)>0, \qquad \forall\xi\in\Vmc\setminus\{0\}. 
\end{equation}
This condition is, in general, stronger than \bla{iv}. It is easy to construct examples 
where \bla{iv} holds and \bla{vii} fails. It is easy to see that the $\mc$-valued scalar product is naturally extendible to higher-order tensor blocks $\Vmc^{\otimes n}$. All 
algebraic properties are preserved. 

In what follows, it will be assumed that conditions \bla{i}--\bla{iii} and \bla{v}--\bla{vii} 
are satisfied. In particular, we see that $\Vmc$ equipped with $\langle\rangle$ gives us 
a (generally unbounded) unitary bimodule over $\mc$. Furthermore, the extended scalar products on spaces $\Vmc^{\otimes n}$ are understandable as $n$-fold tensor iterations of the 
initial bimodule $\Vmc$. In particular, it follows that all extended $\langle\rangle$ are
strictly positive, too (all $\Vmc^{\otimes n}$ are unitary bimodules). 

\begin{lem} The braid operator $\sigma\colon\Vmc^{\otimes 2}\rightarrow\Vmc^{\otimes 2}$
is hermitian. In other words, we have 
\begin{equation}
\langle \psi,\sigma(\xi)\rangle=\langle\sigma(\psi),\xi\rangle\qquad\forall\psi,
\xi\in\Vmc^{\otimes 2}. 
\end{equation}
In particular, the map $\sigma\colon\V\otimes \V\rightarrow\V\otimes\V$ is diagonalizable, 
and has real eigenvalues. 
\end{lem}

\begin{pf} The hermicity property follows from \eqref{*-sigma} and 
\begin{equation}
g( \psi,\sigma[\xi])=g(\sigma[\psi],\xi)\qquad\forall\psi,\xi\in\V\otimes\V
\end{equation}
which, in its turn, follows from the definition of the bimodule structure on $\Vmc$. 
Let $\omega\colon\mc\rightarrow\mathbb{C}$ be an arbitrary faithful state on $\mc$. 
Then $\omega\langle\rangle$ is a scalar product on $\V\otimes\V$, and $\sigma$ is
hermitian with respect to this scalar product. 
\end{pf} 

We are now ready to construct and study braided Clifford algebras. Conceptually, we 
follow \cite{d-ozi} which means that our Clifford algebras will be 
understood as Chevalley-Kahler-type deformations of braided exterior algebras. 

Let $\V^\wedge$ be the braided exterior algebra \cite{w3} built over $(\sigma,\V)$. 
This algebra is defined as $$ \V^\wedge=\V^\otimes/\ker(A_\sigma) $$
where $A_\sigma$ is the total braided antisymmetrizer map. The exterior algebra 
gets its *-structure from $\V^\otimes$. Furthermore, the algebra $\V^\wedge$
possess a natural braided Hopf algebra structure, where the coproduct map is 
specified by 
\begin{equation}
\phi(\alpha)=\sum_{i=0}^n B_{in-i}(\alpha)\qquad \alpha\in \V^{\wedge n}
\end{equation}
and $B_{kl}\colon\V^{\wedge k+l}\rightarrow\V^{\wedge k}\otimes\V^{\wedge l}$
is the corresponding braided inverse-shuffle operator. 

The coproduct map has a particularly elegant form if we make natural
identifications
$$ 
\V^{\wedge n}\leftrightarrow \im(A_\sigma^n)
$$
induced by the antisymmetrizer map. In terms of these identifications, we have
\begin{equation}
\phi(\alpha)=\sum_{i=0}^n \bigl\{x_1\otimes\dots \otimes x_i\bigr\}\otimes
\bigl\{x_{i+1}\otimes\dots \otimes
x_n\bigr\} 
\end{equation}
with $\alpha=\Sum x_1\otimes \dots\otimes x_n$.
The antipode map is braided-antimultiplicative (acting as a total
$\sigma$-inverse permutation). The Hopf algebra structure is compatible
with the *-involution, in the sense that 
$$ \phi(\alpha^*)=\phi(\alpha)^*\qquad \k(\alpha^*)=\k(\alpha)^*. $$

All considerations with the braided exterior algebra are incorporable to 
the level of $\mc$-modules. More precisely, let $\Vmc^\wedge$ be the braided 
exterior algebra constructed from $\Vmc$ and the extended $\sigma$. We have
the following natural identifications
$$ \Vmc^\wedge\leftrightarrow\V^\wedge\otimes\mc\leftrightarrow\mc\otimes\V^\wedge $$
of right/left $\mc$-modules. 
The coproduct map $\phi$ is naturally (and uniquely) extendible to a $\mc$-bilinear 
map $\phi\colon\Vmc^\wedge\rightarrow\Vmc^\wedge\otimes_{\mc}\Vmc^\wedge$. In a similar
way, it is possible to extend the coinverse and the counit. 

Observe now that the block antisymmetrizers 
$A_\sigma^n\colon\V^{\otimes n}\rightarrow\V^{\otimes n}$ are 
hermitian maps, and commute with the *-structure. Hence the pairing
$g_\wedge\colon\Vmc^{\otimes}\times \Vmc^{\otimes}\rightarrow\mc$ defined by
$$g_\wedge(\psi,\xi)=g(\psi, A_\sigma\xi)\qquad g_\wedge(1,1)=1 $$
is projectable down to a map $g_\wedge\colon\Vmc^\wedge\times\Vmc^\wedge\rightarrow\mc$. 

Our braided Clifford algebra $\ClV$ will be 
identified with $\Vmc^\wedge$ at the level of $\mc$-bimodules. The *-structure will also be the same. However, 
$\ClV$ will be equipped with a {\it new product} defined by  
\begin{equation}
\tilde{m}=m(\id\otimes g_\wedge\otimes \id)(\phi\otimes\phi),
\end{equation}
where $m\colon\Vmc^\wedge\otimes\Vmc^\wedge\rightarrow\Vmc^\wedge$ is the original
product in $\Vmc$. In the above formula, we have adopted the standard {\it Cliffordization 
procedure} from the classical theory \cite{rota}. The first thing to examine is that 
we really obtain a nice *-algebra this way. 

\begin{pro} The product $\tilde{m}$ is associative and $1$ is the unit element. Moreover, 
the *-involution is $\tilde{m}$-antimultiplicative. 
\end{pro}

\begin{pf} The associativity of the product follows from braided-multiplicativity of $\phi$, 
property \eqref{funny} and the following interesting identities 
\begin{align}
g_\wedge(m\otimes\id)&=g_\wedge(\id\otimes g_\wedge\otimes\id)(\id\otimes \phi)\\
g_\wedge(\id\otimes m)&=g_\wedge(\id\otimes g_\wedge\otimes\id)(\phi\otimes \id).
\end{align} 
The fact that $1$ is the $\tilde{m}$-unit follows from 
\begin{equation}
g_\wedge(1,\alpha)=g_\wedge(\alpha,1)=\e(\alpha).
\end{equation}
Finally the $\tilde{m}$-antimultiplicativity of $*$ follows from standard commutation 
relations between $*$ and $m, \phi, g_\wedge$. 
\end{pf}

\begin{defn}
The constructed *-algebra $\ClV$ is called the {\it braided Clifford algebra} associated 
to $\Vmc$, $\sigma$ and $g$. 
\end{defn}

We are now going to talk about possible C*-type norms on $\ClV$. For this to work, 
it will be necessary to introduce a last set of our assumptions,
regarding a more subtle behavior of $\sigma$.

A nice way to get such properties is to postulate the existence of
an appropriate {\it auxiliary braid operator}
$\tau\colon\Vmc^{\otimes 2}\rightarrow\Vmc^{\otimes 2}$, as it will
be discussed near the end of this Appendix. Playing with two
braid operators will also enable us to prove interesting properties
of $\sigma$ and its braided exterior algebra $\Vmc^\wedge$.

We shall proceed without making any extra assumptions
on the existence and properties of $\tau$, however we have
to postulate the positivity of braided antisymmetrizer maps.

\bla{viii} {\bf Positivity of braided antisymmetrizers.} All braided antisymmetrizer maps $A_\sigma^n\colon \Vmc^{\otimes n}\rightarrow \Vmc^{\otimes n}$ are positive operators.

The positivity property is crucial to define a C*-algebraic structure on the Clifford
algebra, because only in this case the scalar product $\langle\rangle_\wedge$ on 
$\Vmc^\wedge$ given by the formula
\begin{equation} 
\langle\alpha,\beta\rangle_\wedge=g_\wedge(\alpha^*,\beta)
\end{equation}
will be strictly positive (giving us a structure of a generally unbounded unitary bimodule 
over $\mc$).

\begin{pro}
Let us consider the counit map $\e\colon\ClV\rightarrow\mc$. It is hermitian, $\mc$-bilinear
and strictly positive. 
\end{pro}

\begin{pf} Hermicity and $\mc$-linearity are obvious (the counit here is just the projection on 
$\mc$). The strict positivity follows from the identity
\begin{equation}\label{e-()}
\e(\alpha^*\beta)=\langle\alpha,\beta\rangle_\wedge. \qed
\end{equation}
\renewcommand{\qed}{}
\end{pf}

When dealing with Hilbert space operators, there is an interesting assumption 
we can add to the list of properties of $\mc$---we can assume that the set of
C*-algebraic norms on $\mc$ distinguishes elements of $\mc$. Not every *-algebra
possesses this property, and many *-algebras do not admit any representation by
bounded operators. However if $\mc$ admits C*-algebraic norms, then they would 
be naturally extendible to $\ClV$.

To see this, we can consider the {\it left regular representation} of $\ClV$ in
the $\mc$-bimodule $\Vmc^\wedge$. According to \eqref{e-()},
this representation is a *-representation
$$ \langle \alpha, T\beta\rangle_\wedge=\langle T^*\alpha, \beta\rangle_\wedge\qquad \forall\alpha,\beta\in\Vmc^\wedge\quad
\forall T\in \ClV$$ and it is easy to see that all the operators $T\in \ClV$ are {\it continuous}. This allows us to construct natural C*-algebra norms on $\ClV$.

********\par
Let us now analyze a couple of interesting special cases when the 
positivity of antisymmetrizers would hold automatically. Assume that a selfadjoint
braid operator $\tau\colon\V^{\otimes 2}\rightarrow\V^{\otimes 2}$
is given, satisfying $*\tau *=\tau^{-1}$, extendible by $\mc$-linearity to 
$\Vmc^{\otimes 2}$, and such that
\begin{equation}\label{sigma-tau}
\im(I-\sigma)=\ker(I+\tau),
\end{equation}
or equivalently
\begin{equation}\label{sigma-tau2}
\im(I+\tau)=\ker(I-\sigma). 
\end{equation}
An immediate consequence is that $\sigma$ and $\tau$ commute. Moreover
\begin{equation}\label{imA-tau}
\im(A_\sigma^n)\subseteq\Bigl\{\text{$\tau$-antisymmetric $n$-tensors}\Bigr\}. 
\end{equation}
This inclusion is a simple consequence of \eqref{sigma-tau} and the fact that we can write
$$A_\sigma^n=[\id^k\otimes (I-\sigma)\otimes\id^{n-k-2}]T_k $$
where $T_k\colon\Vmc^{\otimes n}\rightarrow\Vmc^{\otimes n}$
is the part of the antisymmetrizer sum, containing permutations whose
inverse does not reverse the order of $k+1$ and $k+2$.  

\begin{pro} Let us assume that all $\sigma$-twists act in the same
way on the vectors from the space of $\tau$-antisymmetric
$n$-tensors. Then, this space is $\sigma$-invariant.

If in addition $1\in\Bbb{C}$ is the only positive eigenvalue
of $\sigma$, then all braided $\sigma$-antisymmetrizers will be
positive. 
\end{pro}

\begin{pf} The fact that $\sigma$-twists act in the same way on the $\tau$-antisymmetric vectors means that we can always (trivial for $n>4$ for $n=2,3,4$ it is necessary to use the fact that $\sigma$ and $\tau$ commute) exchange them with $\tau$-twists (acting as $-1$).
Hence, the space of $\tau$-antisymmetric tensors is $\sigma$-invariant. 

The second assumption means that $\sigma\colon\ker(I+\tau)\rightarrow\ker(I+\tau)$ is 
strictly negative, as the space $\ker(I+\tau)$ is spanned by all
negative-eigensubspaces of $\sigma$. 

Therefore, all $\sigma$-antisymmetrizers are positive, and their images 
coincide with $\tau$-antisymmetric spaces. 
\end{pf}

\begin{pro} Let us consider mutually equivalent properties
\begin{equation}\label{tau2-sigma}
\begin{aligned}
(\tau\otimes\id)(\id\otimes\tau)(\sigma\otimes\id)&=
(\id\otimes\sigma)(\tau\otimes\id)(\id\otimes\tau)\\
(\sigma\otimes\id)(\id\otimes\tau)(\tau\otimes\id)&=
(\id\otimes\tau)(\tau\otimes\id)(\id\otimes\sigma).
\end{aligned}
\end{equation}
Above equations transform one to another by the *-conjugation. 
If they hold, and if $\tau$-antisymmetric $n$-tensors are invariant
under $\sigma$-twists, then all $\sigma$-twists act in the same way in
this space.
\end{pro}

\begin{pf} Let us consider the case $n=3$. If $\psi\in\V^{\otimes 3}$ is completely
$\tau$-antisymmetric and if the invariance property holds then \eqref{tau2-sigma} gives
$(\id\otimes\sigma)(\psi)=(\sigma\otimes\id)(\psi)$. 
\end{pf}

Furthermore, let us observe that the following strange equalities are
equivalent
\begin{equation}\label{braid1-sigma-tau}
\begin{aligned}
(\tau\otimes\id)(\id\otimes\sigma)(\tau\otimes\id)&=
(\id\otimes\sigma)(\tau\otimes\id)(\id\otimes\tau)\\
(\tau\otimes\id)(\id\otimes\sigma)(\tau\otimes\id)&=
(\id\otimes\tau)(\tau\otimes\id)(\id\otimes\sigma),
\end{aligned}
\end{equation}
as well as the equalities 
\begin{equation}\label{braid2-sigma-tau}
\begin{aligned}
(\id\otimes\tau)(\sigma\otimes\id)(\id\otimes \tau)&=
(\sigma\otimes\id)(\id\otimes\tau)(\tau\otimes\id)\\
(\id\otimes\tau)(\sigma\otimes\id)(\id\otimes \tau)&=
(\tau\otimes\id)(\id\otimes\tau)(\sigma\otimes\id). 
\end{aligned}
\end{equation}
Actually, the equivalent equalities are just mutually adjoint. The
equivalence also follows from the braid equation for $\tau$ and
the fact that $\sigma$ commutes with $\tau$.

\begin{lem} If the above equalities hold, then
the spaces of fully $\tau$-antisymmetric tensors of order
$n\geq 2$ are invariant under actions of all possible $\sigma$-twists.
\end{lem}

\begin{pf}
The invariance under $\sigma$-twists easily follows from the commutation
property, and equalities \eqref{braid1-sigma-tau}--\eqref{braid2-sigma-tau}. 
\end{pf}

\begin{pro} If the spaces of $\tau$-antisymmetric operators are invariant
under all possible $\sigma$-twists and if the restriction 
$\sigma\colon\ker(I+\tau)\rightarrow \ker(I+\tau)$ is negative, then all 
antisymmetrizer maps $A_\sigma^n$ are positive. 

Furthermore, we have
\begin{equation}\label{imA-tau2}
\im(A_\sigma^n)=\Bigl\{\text{$\tau$-antisymmetric tensors}\Bigr\}. 
\end{equation}
for each $n\geq 2$. 
\end{pro}

\begin{pf} At first, let us observe that braided antisymmetrizers satisfy
the following interesting recursive relations
\begin{equation}\label{A-rec}
A_\sigma^{n+1}=\id\otimes A_\sigma^n-(\id\otimes Y_n)(\sigma\otimes A_\sigma^{n-1})(\id\otimes Y_n^\dagger),
\end{equation}
where $Y_n\colon\Vmc^{\otimes n}\rightarrow\Vmc^{\otimes n}$ is given by
$$
Y_n=-\sum_{k=2}^n(-)^k \pi_{kn,\sigma}
$$
and sum goes over permutations $\pi_{kn}\in S_n$ transposing $1$ and
the blocks $\{2,\dots,k\}$ while acting trivially in $\{k+1, \dots, n\}$. 

Now, using induction on $n$, the negativity assumption for $\sigma:\ker(I+\tau)\rightarrow\ker(I+\tau)$, and
recursive formulas \eqref{A-rec} it follows that restricted antisymmetrizers
$$
A_\sigma^n\colon\Bigl\{\text{$\tau$-antisymmetric $n$-tensors}\Bigr\}
\rightarrow \Bigl\{\text{$\tau$-antisymmetric $n$-tensors}\Bigr\}
$$
are {\it strictly positive, in particular invertible} operators. 
In particular \eqref{imA-tau2} holds and obviously $A_\sigma^n$ are
positive everywhere.
\end{pf}

By the way, it is worth mentioning that 
\begin{equation}
\Bigl\{\text{$\tau$-antisymm $n$-tensors}\Bigr\}^\perp=
\sum_k \Vmc^{\otimes k}\otimes \im(I+\tau)\otimes \Vmc^{n-k-2}=\ker(A_\sigma^n). 
\end{equation}
If \eqref{imA-tau2} holds, it follows that the algebra $\Vmc^\wedge$ is
quadratic (generated by its quadratic relations).

In this case, the Clifford algebra $\ClV$ can be viewed as the algebra built over 
$\Vmc$ together with generating relations
\begin{equation}
\sum_\alpha x_\alpha y_\alpha=\sum_\alpha g(x_\alpha,y_\alpha)
\end{equation}
where $\Sum_\alpha \sigma(x_\alpha\otimes y_\alpha)=
\Sum_\alpha x_\alpha\otimes y_\alpha$.

Well, now we can {\it introduce spinors} simply as vectors of irreducible representations of 
the *-algebra $\ClV$ in a finite-dimensional Hilbert space (or more generally, by bounded operators). Recall that every such a representation is (as generally for C*-algebras) obtained from a pure state $\omega\colon\ClV\rightarrow\mathbb{C}$ via the GNS construction. 

In contrast to the classical theory, the algebra $\ClV$ may be infinite dimensional (the most
interesting situations appear when $\mc$ is infinite-dimensional) and possess non-equivalent
irreducible representations. 

In our main context of framed quantum principal bundles, the operator $\tau$ was coming
from the appropriate bicovariant bimodule. In this context it is natural to assume that
$\mc$ is of a `bicovariant nature' too. Specifically, this requires the existence of 
a right $\cal{A}$-module structure $\circ\colon\mc\otimes\cal{A}\rightarrow\mc$ and a 
right $\cal{A}$-comodule structure $\vmc\colon\mc\rightarrow\mc\otimes\cal{A}$ which is a 
(continuous) unital *-homomorphism and such that
\begin{gather*}
(\alpha\beta){\circ a}=(\alpha{\circ} a^{(1)})(\beta{\circ}a^{(2)})\qquad 1{\circ}a =\e(a) 1\\
\mc(q{\circ} a)=\sum_\alpha (q_\alpha{\circ} a^{(2)})\otimes \k(a^{(1)})c_\alpha a^{(3)}
\end{gather*}
The maps $\circ$ and $\vmc$ are completely fixed by postulating
\begin{gather*}
g(x,y)\circ a=g(x{\circ} a^{(1)},y{\circ} a^{(2)})\\
\vmc g(x,y)=(g\otimes\id)\v(x\otimes y). 
\end{gather*}
The above formulas extend to the elements from $\V^\wedge$ straightforwardly. 

\section{Intertwiner Bimodules $\&$ Vector Bundles}

The aim of this appendix is to sketch basic ideas and constructions related to associated vector 
bundles. We shall mainly follow the exposition of \cite{d-tann}. 

Let $\RepG$ be the category of finite-dimensional representations of a compact quantum group 
$G$. The objects of $\RepG$ are finite-dimensional representations of $G$ and the arrows are 
the intertwiners between the corresponding representations. For each $u\in\RepG$ let us denote 
by $H_u$ the corresponding carrier vector space. 

Let us consider an arbitrary quantum principal $G$-bundle $P=(\cal{B},i,F)$ over a quantum space 
$M$. Let us denote by $\bim{u}=\Mor(u,F)$ be the space of intertwiners between a 
given representation $u\colon H_u\rightarrow H_u\otimes\cal{A}$ 
and $F\colon\cal{B}\rightarrow\cal{B}\otimes\cal{A}$. 
The spaces $\bim{u}$ are $\cal{V}$-bimodules, in a natural manner. We have 
$\bim{\varnothing}\leftrightarrow\cal{V}$. 

In accordance with classical geometry, the spaces $\bim{u}$ are interpretable as (the smooth 
sections of) the {\it associated vector bundles}.  It is possible to give an important 
alternative interpretation of $\bim{u}$ as certain invariant subspaces. More precisely, 
let us consider the contragradient representation $u^c\colon H_u^*\rightarrow H^*_u\otimes\cal{A}$. 
Then we have a natural identification
\begin{equation}
\bim{u}\leftrightarrow\Bigl\{\psi\in\cal{B}\otimes H_u^*\Bigm\vert (F{\times}u^c)(\psi)=\psi\otimes 1\Bigr\}.
\end{equation}
Explicitly, the identification is given by
\begin{equation}\label{interpret-bim}
\begin{gathered}
\bim{u}\ni\psi\longmapsto \sum_i\psi(e_i)\otimes e_i^*\\
\sum_\alpha b_\alpha\otimes f_\alpha\longmapsto \Bigl\{x\mapsto 
\sum_\alpha b_\alpha f_\alpha(x)\Bigr\}, 
\end{gathered}
\end{equation}
where $e_i$ form a basis (say, orthonormal) in $H_u$ and $e_i^*\in H_u^*$ are 
the corresponding biorthogonal vectors. 
In various considerations it comes very handy to pass from one interpretation to another. 

The following natural bimodule isomorphism holds:
\begin{equation}
\bim{u\times v}\leftrightarrow\bim{u}\otimes_{\cal{V}}\!\bim{v} \qquad\quad \forall u,v\in\RepG.
\end{equation}
This isomorphism is induced by the product in $\cal{B}$, via
$\varphi\otimes\psi\colon x\otimes y\mapsto\varphi(x)\psi(y)$.  
Furthermore, every map $f\in\Mor(u,v)$ induces, via the composition of intertwiners, 
a bimodule homomorphism $\map{f}\colon\bim{v}\rightarrow\bim{u}$. 

We have a system of natural bimodule anti-isomorphisms $\amaph{u}\colon
\bim{u}\rightarrow\bim{\bar{u}}$, defined by the diagram
\begin{equation}
\begin{CD}
H_{u} @>{\mbox{$\phantom{\amaph{u}}\varphi$}}>> \cal{B}\\
@V{\mbox{$j_u$}}VV @VV{\mbox{$*$}}V\\
H_{u}^* @>>{\mbox{$\amaph{u}\varphi$}}> \cal{B}
\end{CD}
\end{equation}
In such a way we can interpret $\bim{\bar{u}}$ as the {\it conjugate 
bimodule} of $\bim{u}$.

Let us now focus on the main context of this paper, and consider a graded *-algebra
$\hor_P$ playing the role of abstract `horizontal forms'. Here we shall assume that
$\hor_P^0=\cal{B}$ and that there exists a coassociative counital *-homomorphism 
$\Fh\colon\hor_P\rightarrow\hor_P\otimes\cal{A}$ extending the map $F$. Let $\Omega_M$ be
the corresponding $\Fh$-fixed point subalgebra. 

Applying similar intertwiner considerations to the algebra $\hor_P$
leads to $\Omega_M$-bimodules $\hbim{u}=\Mor(u,F^\wedge)$. 
These spaces are naturally graded, and obviously
$\hbim{u}^0=\bim{u}$. Moreover, it can be shown that canonically 
\begin{equation}
\bim{u}\otimes_{\cal{V}}\!\Omega_M\leftrightarrow\hbim{u}\leftrightarrow
\Omega_M\otimes_{\cal{V}}\!\bim{u}.
\end{equation}
These decompositions are induced by the bimodule product in $\hbim{u}$. 
The spaces $\hbim{u}$ are alternatively viewable in the same way \eqref{interpret-bim}. 

The structure of $\hbim{u}$ is expressible in terms of
$\bim{u}$. Composing the above two identifications we obtain canonical
flip-over maps $$\sigma_u\colon\bim{u}\otimes_{\cal{V}}\!\Omega_M
\rightarrow\Omega_M\otimes_{\cal{V}}\!\bim{u}.$$ These maps are
grade-preserving and act as identity on $\bim{u}$. The
$\Omega_M$-bimodule structure on $\hbim{u}$ is expressed by the
following equation
\begin{equation}
\sigma_u(\id\otimes m_*)=(m_*\otimes\id)
(\id\otimes\sigma_u)(\sigma_u\otimes\id),
\end{equation}
where $m_*$ is the product in $\Omega_M$. Intertwiner homomorphisms
$\map{f}$ and conjugation maps $\amaph{u}$ obey the following diagrams
\begin{equation*}
\begin{CD} \bim{v}\otimes_{\cal{V}}\!\Omega_M @>{\mbox{$\sigma_v$}}>>
\Omega_M\otimes_{\cal{V}}\!\bim{v}\\
@V{\mbox{$\map{f}\otimes\id$}}VV @VV{\mbox{$\id\otimes\map{f}$}}V\\
\bim{u}\otimes_{\cal{V}}\!\Omega_M @>>{\mbox{$\sigma_u$}}>
\Omega_M\otimes_{\cal{V}}\!\bim{u}
\end{CD}\qquad
\begin{CD} \bim{u}\otimes_{\cal{V}}\!\Omega_M @>{\mbox{$\sigma_u$}}>>
\Omega_M\otimes_{\cal{V}}\!\bim{u}\\
@V{\mbox{$\amaph{u}$}}VV  @VV{\mbox{$\amaph{u}$}}V\\
\Omega_M\otimes_{\cal{V}}\!\bim{\bar{u}} @<<{\mbox{$\sigma_{\bar{u}}$}}<
\bim{\bar{u}}\otimes_{\cal{V}}\!\Omega_M
\end{CD}
\end{equation*}
In fact, the first diagram characterizes elements of
$\hom(\bim{v},\bim{u})$ that are extendible to corresponding
$\Omega_M$-bimodule
homomorphisms (the left and right $\Omega_M$-linear extensions coincide).

We have the following natural decomposition of the algebra of horizontal forms 
\begin{equation}
\begin{gathered}
\hor_P=\sideset{}{^\oplus}\sum_{\alpha\in\cal{T}}\cal{H}^\alpha(P)\qquad
\cal{H}^\alpha(P)=\hbim{\alpha}\otimes H_\alpha\\
\hor_P\leftrightarrow\Omega_M\otimes_{\cal{V}}\!\cal{B}
\leftrightarrow\cal{B}\otimes_{\cal{V}}\!\Omega_M.
\end{gathered}
\end{equation}

Let us assume that $\hor_P$ is equipped with an $\Fh$-invariant scalar product {$<>$}, ensuring
that the spaces $\hbim{\alpha}$ are mutually orthogonal. Then it is possible to naturally induce
scalar products in all intertwiner $\Omega_M$-bimodules. The induced product is 
given by 
\begin{equation}
<\varphi,\psi>=\sum_{ij} [C_u]_{ji}<\varphi(e_i),\psi(e_j)>
\end{equation} 
where $e_i$ form an orthonormal basis in $H_u$. In the alternative picture, the scalar product is simply
given by tensoring the scalar products in $\hor_P$ and $H_u^*$. In the case of quantum 
Riemannian/spin manifolds, the above scalar product reads
\begin{equation}
<\varphi,\psi>=\int_M\Bigl\{\sum_{ij}[C_u]_{ji}\varphi(e_i)^*\hstrP[\psi(e_j)]\Bigr\}. 
\end{equation} 
Let us suppose that we have a linear operator $T\colon\hor_P\rightarrow\hor_P$
which intertwines the right action $\Fh$. Then it acts naturally in all intertwiner 
bimodules $\hbim{u}$. The action is defined by simply taking the composition with the
intertwiners. In other words
$$ [T_u\psi](x)=T[\psi(x)] \qquad\quad \psi\in\hbim{u}, x\in H_u$$
or in the alternative picture 
$$ T_u\Bigl\{\sum_\alpha\varphi_\alpha\otimes f_\alpha\Bigr\}=\Sum_\alpha T(\varphi_\alpha)
\otimes x_\alpha. $$

It is clear that in such a way all algebraic relations between $\Fh$-covariant operators in
$\hor_P$ are preserved. Furthermore, the adjoint operation is preserved, in a natural manner. 
More precisely, if $T$ is formally adjointable then $T_u$ is formally adjointable too, 
for each $u\in\RepG$ and $[T_u]^\dagger=[T^\dagger]_u$. 
This is a direct consequence of the definition of the scalar product in our bimodules.

\end{document}